\documentclass[11pt]{article}
\usepackage{amsmath,amsthm,amsfonts,amssymb,tikz-cd}
\usepackage[colorlinks,linktoc=page]{hyperref}

\title{A categorical $\mathfrak{sl}_2$ action \\ on some moduli spaces of sheaves}
\author{Nicolas Addington and Ryan Takahashi}
\date{}

\DeclareMathOperator \ST {ST} 
\DeclareMathOperator \Pic {Pic}
\DeclareMathOperator \Gr {Gr}
\DeclareMathOperator \supp {supp}
\DeclareMathOperator \Hom {Hom}
\DeclareMathOperator \codim {codim}
\DeclareMathOperator \rank {rank}
\DeclareMathOperator \Ext {Ext}
\DeclareMathOperator \sheafExt {\mathcal E xt}
\DeclareMathOperator \cone {cone}
\DeclareMathOperator \Hilb {Hilb}
\DeclareMathOperator \Br {Br}
\DeclareMathOperator \RGamma {R\Gamma}
\DeclareMathOperator \Spec {Spec}
\DeclareMathOperator \relProj {\underline{Proj}}
\DeclareMathOperator \Sym {Sym}

\DeclareMathOperator \sheafHom {\mathcal H om}
\DeclareMathOperator \Bl {Bl}

\newcommand \C {\mathbb C}
\renewcommand \H {\mathcal H}
\renewcommand \O {\mathcal O}
\renewcommand \P {\mathbb P}
\newcommand \Z {\mathbb Z}

\newcommand \TGr {\mathrm{T^*Gr}}

\renewcommand \sl {\mathfrak{sl}}
\renewcommand \th {\text{th}}

\newtheorem*{mainthm*}{Main Theorem}
\newtheorem*{maincor*}{Corollary}

\newtheorem{thm}{Theorem}[section]
\newtheorem{prop}[thm]{Proposition}
\newtheorem{lem}[thm]{Lemma}
\newtheorem{cond}{Condition}

\theoremstyle{definition}
\newtheorem{defn}[thm]{Definition}
\newtheorem{rmk}[thm]{Remark}

\numberwithin{equation}{section}

\begin{document}
\maketitle

\begin{abstract}
We study certain sequences of moduli spaces of sheaves on K3 surfaces, building on work of Markman, Yoshioka, and Nakajima.  We show that these sequences can be given the structure of a geometric categorical $\sl_2$ action in the sense of Cautis, Kamnitzer, and Licata.  As a corollary, we get an equivalence between derived categories of some moduli spaces that are birational via stratified Mukai flops.
\end{abstract}

\section*{Introduction}

\setlength \parskip \smallskipamount

Equivalences between derived categories of K3 surfaces induce birational maps between moduli spaces of stable sheaves, with descriptions ranging from straightforward to very complicated.  \emph{Stratified Mukai flops}, introduced by Markman \cite{markman_brill_noether} and Yoshioka \cite{yoshioka_reflection}, are a rich but still manageable class that arise as follows.  Let $S$ be a complex K3 surface with $\Pic(S) \cong \Z$ generated by an ample line bundle $h$.  Let $\ST_\O \colon D(S) \to D(S)$ be the spherical twist around $\O_S$, introduced by Mukai \cite[\S2]{mukai_tata} under the name ``reflection functor'' and since studied and generalized by many authors.  Then $\ST_\O$ induces a birational map
\begin{equation} \label{the_flop} \tag{$*$}
\begin{tikzcd}[row sep = tiny]
M(r,h,s)\footnotemark \arrow[leftrightarrow,dashed]{r} &M(-s,h,-r) &\text{if $s \le 0$,} \\
M(r,h,s) \arrow[leftrightarrow,dashed]{r} &M(s,h,r) &\text{if $s \ge 0$.}
\end{tikzcd}
\end{equation}
\footnotetext{The triples in parentheses are Mukai vectors, so $M(r,h,s)$ parametrizes sheaves of rank $r$, first Chern class $c_1 = h$, and Euler characteristic $\chi = r+s$.}%
These moduli spaces are smooth, projective, holomorphic symplectic varieties.  The map is biregular on the open set parametrizing sheaves whose cohomology is as small as possible, and indeterminate on the \emph{Brill--Noether locus} parametrizing sheaves whose cohomology jumps up: $\ST_\O$ 
might take these to unstable sheaves, 
or to 2-term complexes of sheaves. 
In the simplest cases, the Brill--Noether locus is a Lagrangian $\P^n$ or coisotropic $\P^n$ bundle, and the birational map is an ordinary Mukai flop, which blows up the $\P^n$ (bundle) and then blows down in another way \cite[\S3]{mukai_inventiones}.  But in general the birational map has a more complicated structure, which Markman analyzed in detail.  Yoshioka allowed spherical twists around other rigid stable vector bundles, but for clarity we restrict our attention to $\O_S$.

The local model for a stratified Mukai flop is a birational map between cotangent bundles of dual Grassmannians
\[ \begin{tikzcd}
\TGr(k,n) \arrow[leftrightarrow,dashed]{r} & \TGr(n-k,n),
\end{tikzcd} \]
which are non-compact symplectic varieties.  In fact they are Nakajima quiver varieties, and Nakajima showed that $\bigoplus_k H^*(\TGr(k,n), \C)$ carries an action of the affine Lie algebra $\widehat\sl_2$.  In \cite{nakajima_k3}, he constructed an analogous $\sl_2$ action on the cohomology groups of the compact moduli spaces discussed above, drawing on Markman and Yoshioka's work, and conjectured that these actions should lift to derived categories of coherent sheaves.

In \cite{ckl1795} and \cite{ckl1796}, Cautis, Kamnitzer, and Licata introduced the notion of a \emph{geometric categorical $\sl_2$ action}, and for their main example constructed such an action on the derived categories of coherent sheaves on the sequence of spaces
\begin{align*}
\dots && \TGr(k-1,n) && \TGr(k,n) && \TGr(k+1,n) && \dots.
\end{align*}
As a corollary, they obtained in \cite{ckl1797} an interesting equivalence
\[ D(\TGr(k,n)) \cong D(\TGr(n-k,n)), \]
generalizing the equivalence for ordinary Mukai flops (the case $k=1$) due to Namikawa \cite{namikawa} and Kawamata \cite[\S5]{kawamata}.

In this paper we construct such an action on the compact moduli spaces discussed above.

\begin{mainthm*}
Let $S$ be a complex K3 surface with $\Pic(S) \cong \Z$ generated by an ample line bundle $h$.  Fix integers $r$ and $s$, and consider the moduli spaces of stable sheaves on $S$
\begin{align*}
\dots && M(r-1,h,s-1) && M(r,h,s) && M(r+1,h,s+1) && \dots,
\end{align*}
with a correction described in \S\ref{sectioncorrespondences} when the rank becomes negative on the left.

There are natural correspondences over each pair of moduli spaces, and line bundles on the correspondences, that give a collection of Fourier--Mukai functors with the structure of a geometric categorical $\sl_2$ action.
\end{mainthm*}
\pagebreak 

Only finitely many of the moduli spaces in question are non-empty.  We do not require that they admit universal sheaves, which would translate to the numerical condition $\gcd(r-s, h^2) = 1$.  Our correspondences are the same ones that Nakajima used in \cite{nakajima_k3}, but the need to choose line bundles on them is new in the derived setting, and quite delicate.

Cautis, Kamnitzer, and Licata's definition sidesteps the need to check a complicated set of ``nil affine Hecke relations'' by choosing a 1-parameter deformation of each space, such that the rest of the $\sl_2$ structure does not deform along with the spaces.  In their example, the deformation comes from the unique non-split extension of $\O_{\Gr(k,n)}$ by $\TGr(k,n)$.  In ours, the deformation is the twistor family of a hyperk\"ahler metric on $M$.

As a corollary to our main theorem, we obtain a case of Bondal, Orlov, and Kawamata's conjecture on flops and derived categories:
\begin{maincor*}
The birational moduli spaces appearing in \eqref{the_flop} have equivalent derived categories of coherent sheaves.
\end{maincor*}
\noindent Halpern-Leistner \cite{dhl_magic_windows} has shown that \emph{any} pair of birational moduli spaces of sheaves on K3 surfaces are derived equivalent, using very different methods.  It would be interesting to study any differences between his equivalence and ours, as the groupoid of derived equivalences between varieties related by a flop is often rich in both geometry and representation theory.

Another related work is Jiang and Leung's paper \cite{jl} on derived categories of projectivizations of sheaves that are not vector bundles but have projective dimension 1, with applications to equivalences between varieties related by a flop, and ``flop-flop = twist'' results.  Some of our correspondences are projectivizations of the kind studied by Jiang and Leung.  In general, our correspondences are Grassmannians of (twisted) sheaves of projective dimension 1, and a study of derived cartegories of such Grassmannians along the lines of Kapranov's classic paper \cite{kapranov_gr} might yield further interesting ``flop-flop = twist'' results for our moduli spaces.

Yoshioka \cite{yoshioka_ade} generalized Nakajima's $\sl_2$ action \cite{nakajima_k3} to an action of higher-rank Lie algebras by considering not just one spherical vector bundle, but several of them in an ADE configuration.  In \cite{ckl_quiver}, Cautis, Kamnitzer, and Licata categorified the Lie algebra action on Nakajima quiver varieties of simply-laced type.  It would be natural to try to categorify Yoshioka's action using ideas from \cite{ckl_quiver} and this paper.

Negu\c{t} \cite{negut} categorified Baranovsky's Heisenberg action on cohomology of moduli spaces of higher-rank sheaves on K3, Abelian, or del Pezzo surfaces \cite{baranovsky}, which in turn generalized Nakajima and Grojnowski's action on cohomology of Hilbert schemes of points on surfaces \cite{nakajima, grojnowski}.  It would be interesting to explore interactions between Negu\c{t}'s action and ours.

In \S\S\ref{sectionsl2}--\ref{sectionbn} we review the definition of a geometric categorical $\sl_2$ action and some needed facts about moduli spaces of sheaves on K3 surfaces.  In \S\S\ref{sectioncorrespondences}--\ref{sectionlinebundles} we construct the correspondences between our moduli spaces, give an example to illustrate the construction, and construct the line bundles on the correspondences.  In \S\S\ref{sectioneasycond}--\ref{sectioncond4} we show that our constructions satisfy the required conditions.  In Appendix \ref{appendix} we collect some results on Grassmannians of coherent sheaves.

\subsection*{Acknowledgements} We thank Will Donovan and Ben Elias for helpful discussions, Eyal Markman for advice on \cite{markman_brill_noether}, S\'andor Kov\'acs for advice on rational singularities, Hiraku Nakajima and Andrei Negu\c{t} for bringing \cite{nakajima_k3} and \cite{yoshioka_ade} to our attention, and the referees for their helpful comments.  Both authors were partly supported by NSF grant no.\ DMS-1902213.

\subsection*{Conventions}  $D(X)$ denotes the bounded derived category of coherent sheaves on $X$.  All functors are implicitly derived.  In particular, we write $f_*$ rather than $Rf_*$, and when we need the underived pushforward we write $R^0 f_*$.

Wherever the cohomology of $\P^n$ or a Grassmannian appears, its grading is centered around zero: thus $H^*(\P^1) = \C[1] \oplus \C[-1]$ rather than $\C \oplus \C[-2]$.

A $\P^n$ bundle or Grassmannian bundle is locally trivial in the analytic or \'etale topology, but not necessarily in the Zariski topology.  Other authors might prefer to call this a $\P^n$ fibration.

For $k \le n$, we consider both the Grassmannian $\Gr(k,n)$ of $k$-dimensional subspaces of $\mathbb C^n$, and the dual Grassmannian $\Gr(n,k)$ of $k$-dimensional quotients of $\mathbb C^n$.  These are related via $\Gr(n,k) = \Gr(n-k,n)$.  It will be clear from context which of the two numbers is greater.

\setlength \parskip {0pt}

\pagebreak 

\tableofcontents

\section{Geometric categorical \texorpdfstring{$\sl_2$}{sl\_2} actions}
\label{sectionsl2}

We review Cautis, Kamnitzer, and Licata's definition of geometric categorical $\sl_2$ action, adjusting the notation slightly and assuming that the base field is $\C$. Following \cite[Rmk.~2.6]{ckl1796}, we dispense with the $\C^\times$ action, which is important in their non-compact example but superfluous in our example.

\begin{defn}[{\cite[Def.~2.2.2]{ckl1796}}]
\label{maindef}
A \emph{geometric categorical $\sl_2$ action} consists of the following data.
\begin{enumerate}
\renewcommand \labelenumi {(\roman{enumi})}
\item A sequence of smooth varieties $M_{-N}, M_{-N+1}, \dots, M_{N-1}, M_N$ over $\C.$
\item Fourier-Mukai kernels
\[E^{(k)}_\chi \in D(M_{\chi-k} \times M_{\chi+k}) \text{ and } F^{(k)}_\chi \in D(M_{\chi+k} \times M_{\chi-k}).\]
We write $E_\chi$ for $E^{(1)}_\chi$ and take $E^{(0)}_\chi = \O_\Delta$, and similarly with $F$.
\item For each $M_\chi$, a flat deformation $\widetilde M_\chi \to \C.$
\end{enumerate}
These data are required to satisfy the following conditions.
\begin{enumerate}
\renewcommand \labelenumi {(\roman{enumi})}
\renewcommand \theenumi \labelenumi

\item \label{condi} Each Hom space between two objects of $D(M_\chi)$ is finite dimensional.

\item \label{condii} $E^{(k)}_\chi$ and $F^{(k)}_\chi$ are left and right adjoints of each other up to shift. More precisely,
\begin{enumerate}
\item The right adjoint of $E^{(k)}_\chi$ is $F^{(k)}_\chi[k\chi]$
\item The left adjoint of $E^{(k)}_\chi$ is $F^{(k)}_\chi[-k\chi]$.
\end{enumerate}

\item \label{condiii} At the level of cohomology of complexes we have
\[\H^*(E_{\chi+k} \ast E^{(k)}_{\chi-1}) \cong E^{(k+1)}_\chi \otimes_\C H^*(\P^r),\]
where the grading of $H^*(\P^r)$ is centered around 0.

\item \label{condiv} If $\chi \leq 0$ then \[F_{\chi+1} \ast E_{\chi+1} \cong E_{\chi-1} \ast F_{\chi-1} \oplus P,\] where $\H^*(P) \cong \O_\Delta \otimes_\C H^*(\P^{-\chi-1}).$

Similarly, if $\chi \geq 0$ then \[E_{\chi-1} \ast F_{\chi-1} \cong F_{\chi+1} \ast E_{\chi+1} \oplus P',\] where $\H^*(P') \cong \O_\Delta \otimes_\C H^*(\P^{\chi-1}).$

\item \label{condv} We have
\[\H^*(i_{23*} E_{\chi+1} \ast i_{12*} E_{\chi-1}) \cong E^{(2)}_\chi[-1] \oplus E^{(2)}_\chi[2],\]
where $i_{12}$ and $i_{23}$ are the closed embeddings
\[i_{12}: M_{\chi-2} \times M_\chi \to M_{\chi-2} \times \widetilde M_\chi \]
\[i_{23}: M_\chi \times M_{\chi+2} \to \widetilde M_\chi \times M_{\chi+2}.\]

\item \label{condvi} If $\chi \leq 0$ then for $k' > k$, the image of $\supp(E^{(k)}_{\chi-k})$ under the projection to $M_\chi$ is not contained in the image of $\supp(E^{(k')}_{\chi-k'})$ also under the projection to $M_\chi.$ Similarly, if $\chi \geq 0$ then for $k' > k$, the image of $\supp(E^{(k)}_{\chi+k})$ in $M_\chi$ is not contained in the image of $\supp(E^{(k')}_{\chi+k'})$.

\item \label{condvii} All $E^{(k)}$s and $F^{(k)}$s are sheaves (that is, complexes supported in degree zero).
\end{enumerate}
\end{defn}

It will be helpful to have the following reformulation of condition \ref{condiv}, which Cautis, Kamnitzer, and Licata also use in \cite[\S6.3, step 4]{ckl1795}:
\begin{lem} \label{weak_fours}
To verify condition \ref{condiv}, it is enough to produce exact triangles
\begin{align*}
P &\to F_{\chi+1} * E_{\chi+1} \to E_{\chi-1} * F_{\chi-1} & \text{for }\chi \le 0, \\
P' &\to E_{\chi-1} * F_{\chi-1} \to F_{\chi+1} * E_{\chi+1} & \text{for }\chi \ge 0.
\end{align*}
in $D(M_\chi \times M_\chi)$, because any such triangle is necessarily split.
\end{lem}

\begin{proof} We prove it for $\chi \ge 0$; the other case is similar.  From condition \ref{condii} we find that the right adjoint of $F_{\chi+1}$ is $E_{\chi+1}[-\chi-1]$, so the set of extensions
\[ \Hom_{M_\chi \times M_\chi}(F_{\chi+1} * E_{\chi+1},\ P'[1]) \]
can be rewritten as
\begin{equation} \label{ext1}
\Hom_{M_\chi \times M_{\chi+2}}(E_{\chi+1},\ E_{\chi+1} * P'[-\chi]). 
\end{equation}
If $P'$ is formal, that is, if
\[ P' = \O_\Delta[-\chi+1] \oplus \O_\Delta[-\chi+3] \oplus \dotsb \oplus \O_\Delta[\chi-1], \]
then \eqref{ext1} is
\begin{equation} \label{vanishing_terms}
\Hom(E_{\chi+1},\ E_{\chi+1}[-2\chi+1]) \oplus \dotsb \oplus \Hom(E_{\chi+1},\ E_{\chi+1}[-1]),
\end{equation}
which vanishes because $E_{\chi+1}$ is a sheaf (condition \ref{condvi}).  If $P'$ is not formal then there is a Grothendieck spectral sequence
\[ E_2^{p,q} = \Hom(F_{\chi+1} * E_{\chi+1}, \H^q(P')[p]) \Longrightarrow \Hom(F_{\chi+1} * E_{\chi+1}, P'[p+q]), \]
and the terms appearing on the diagonal $p+q=1$ are the summands of \eqref{vanishing_terms}, which still vanish, so nothing can contribute to \eqref{ext1}.
\end{proof}

The Corollary in the introduction follows from Cautis, Kamnitzer, and Licata's results \cite[Thm.~2.5]{ckl1796} and \cite[Thm.~2.8]{ckl1797}, which together say that a geometric categorical $\sl_2$ action gives rise to an equivalence $D(M_\chi) \cong D(M_{-\chi})$ for each $\chi \ge 0$.  The kernel that induces the equivalence is obtained as a convolution
\[ \dotsb \to F^{(\chi+2)}_{-2} * E^{(2)}_{\chi+2}[-2]
\to F^{(\chi+1)}_{-1} * E^{(1)}_{\chi+1} [-1]
\to F^{(\chi)}_0. \]

\section{Brill--Noether stratification of moduli spaces}
\label{sectionbn}

Fix, once and for all, a complex K3 surface $S$ with $\Pic(S) \cong \Z$ generated by an ample line bundle $h$.  We review a few facts about moduli spaces of stable sheaves on $S$, and refer to \cite[Ch.~10]{huybrechts_k3} for a textbook account.

Let $M(r,h,s)$ denote the moduli space of Gieseker $h$-stable sheaves with Mukai vector
\[ (r,h,s) \in H^0(S,\Z) \oplus H^2(S,\Z) \oplus H^4(S,\Z), \]
where the first component is the rank, the second component is the first Chern class, and the third component $s$ is $\frac12 c_1^2 - c_2 + r$.  The holomorphic Euler characteristic of such a sheaf is
\[ \chi = r+s. \]
Our assumption on the first Chern class ensures that every semi-stable sheaf is stable; thus the moduli space is a smooth, projective, irreducible holomorphic symplectic variety of dimension
\[ \dim M(r,h,s) = h^2 - 2rs + 2 \]
if that number is non-negative, or is empty otherwise.

The \emph{Brill--Noether stratification} of $M(r,h,s)$ is given by the jumping of the cohomology of the sheaves in question.  Observe that if $F$ is a stable sheaf with Mukai vector $(r,h,s)$ then $H^2(F) \cong \Hom(F,\O)^*$ vanishes: if $r = 0$ then $F$ is torsion while $\O$ is torsion-free, and if $r > 0$ then the slope of $F$ is positive while the slope of $\O$ is zero.\footnote{We thank Emanuele Macr\`i for this easy but crucial remark.}  But $H^0(F)$ and $H^1(F)$ may jump as $F$ varies.

If $\chi \ge 0$ then a general point of $M(r,h,s)$ represents a sheaf $F$ with $h^0(F) = \chi$ and $h^1(F) = 0$, and we set
\[ _tM = \{ F \in M(r,h,s) : h^0(F) \ge \chi + t \text{ and } h^1(F) \ge t \}. \]
If $\chi \le 0$ then a general point of $M(r,h,s)$ represents a sheaf $F$ with $h^0(F) = 0$ and $h^1(F) = \left|\chi\right|$, and we set
\[ _tM = \{ F \in M(r,h,s) : h^0(F) \ge t \text{ and } h^1(F) \ge \left|\chi\right| + t \}. \]
Markman \cite[Cor.~34]{markman_brill_noether} and Yoshioka \cite[Lem.~2.6 and Rmk.~2.3]{yoshioka_reflection} showed that $_tM$ is an integral scheme
with
\[ \codim (_tM) = t(\left|\chi\right| + t), \]
if that number is less than or equal to \emph{half} of $\dim M$, and is empty otherwise.


\section{Correspondences}
\label{sectioncorrespondences}

\subsection*{Avoiding negative ranks}

Fix an integer $x$.  Suppose to begin with that $x > h^2/2$, and consider the sequence of moduli spaces
\begin{align*}
M(0,h,-x) && M(1,h,1-x) && M(2,h,2-x) && \dots && M(x,h,0).
\end{align*}
We cannot continue the sequence any further to the right, because the moduli spaces would become empty.  In Definition \ref{maindef}, we set
\[ M_\chi := \begin{cases}
M(r,\,h,\,r-x) & \text{if $\chi = 2r-x$ with $0 \le r \le x$,} \\
\varnothing & \text{otherwise}.
\end{cases} \]

For any integer $\chi$ and positive integer $k$ such that $M_{\chi-k}$ and $M_{\chi+k}$ are not empty, we construct a correspondence
\[ \begin{tikzcd}
& X \ar[swap]{dl}{f} \ar{dr}{g} & \\ M_{\chi-k} & & M_{\chi+k}.
\end{tikzcd} \]
We deliberately refrain from introducing a heavy notation for the correspondences $X$, because we will only ever consider one correspondence between any given pair of moduli spaces. \bigskip

Somewhat informally, we let $X$ be the space of pairs consisting of a sheaf $F \in M_{\chi-k}$ and a quotient $H^1(F) \twoheadrightarrow \C^k$, or equivalently, a sheaf $G \in M_{\chi+k}$ and a subspace $\C^k \subset H^0(G)$.  Let us explain why these are equivalent.  We necessarily have $k \le \rank(G)$, so by \cite[Lem.~25(3)]{markman_brill_noether} or \cite[Lem.~2.1(2-1)]{yoshioka_reflection}, a subspace $\C^k \subset H^0(G)$ determines an injection $\O_S^k \hookrightarrow G$, and the cokernel $F$ is stable.  We have $\chi(F) = \chi-k$ because $\chi(\O_S) = 2$.  From the short exact sequence
\begin{equation} \label{ogf}
\begin{tikzcd}
0 \ar{r} & \O_S^k \ar{r} & G \ar{r} & F \ar{r} & 0
\end{tikzcd}
\end{equation}
we get a long exact sequence in cohomology
\[ \begin{tikzcd}
0 \ar[r] & \C^k \ar[r] & H^0(G) \ar[r] & H^0(F) \ar[lld, controls={+(3,-0.5)and +(-3,.5)}] \\
& 0 \ar[r] & H^1(G) \ar[r] & H^1(F) \ar[lld, controls={+(3,-0.5)and +(-3,.5)}] \\
& \C^k \ar[r] & 0 \ar[r] & 0 \ar[r] & 0,
\end{tikzcd} \]
so $F$ comes with a surjection $H^1(F) \twoheadrightarrow \C^k$.  The opposite direction is much the same: a quotient $H^1(F) \twoheadrightarrow \C^k$ is dual to a subspace $\C^k \subset \Ext^1(F,\O_S),$ which gives an extension of $F$ by $\O_S^k$, that is, the same short exact sequence \eqref{ogf}.  The sheaf $G$ comes with $k$ sections, and is stable by \cite[Lem.~25(1)]{markman_brill_noether} or \cite[Cor.~2.2]{yoshioka_reflection}. \bigskip

More formally, let $U_{\chi-k}$ be a universal sheaf on $S \times M_{\chi-k}$, possibly twisted by a Brauer class pulled back from $M_{\chi-k}$, and let $q$ be the projection from $S \times M_{\chi-k}$ onto $M_{\chi-k}$.  Similarly, let $U_{\chi+k}$ be a (possibly twisted) universal sheaf on $S \times M_{\chi+k}$, and let $q'$ be the projection onto $M_{\chi-k}$.  Then we let $X$ be on the one hand the Grassmannian of $k$-dimensional quotients
\[ X := \Gr(R^1 q_* U_{\chi-k},\,k), \]
with its natural map to $M_{\chi-k}$, and on the other hand
\[ X := \Gr(\sheafExt^2_{q'}(U_{\chi+k},\O),\,k), \]
with its natural map to $M_{\chi+k}$.  We collect some background and needed facts about Grassmannians of coherent sheaves in Appendix \ref{appendix}; it is helpful to know that the fiber of the Grassmannian of a sheaf is the Grassmannian of the fiber of the sheaf.  In our case, the fiber of $R^1 q_* U_{\chi-k}$ over $F \in M_{\chi-k}$ really is $H^1(F)$, because $H^2(F)$ vanishes, although the same cannot be said about $R^0 q_* U_{\chi-k}$: see \cite[III Thm.~12.11]{hartshorne}.  Similarly, the fiber of $\sheafExt^2_{q'}(U_{\chi+k},\O)$ over $G \in M_{\chi+k}$ is really $\Ext^2(G,\O) \cong H^0(G)^*$, although the same cannot be about $\sheafExt^1_{q'}$ or $\sheafExt^0_{q'}$.

The claim, then, is that these two Grassmannians are isomorphic.  Markman proves this in \cite[Thm.~33]{markman_brill_noether},\footnote{For the reader who wants to look closely at this reference, we remark that Markman's $G^0$ and $G_1$ are defined toward the beginning of [ibid., \S5.4].} or see \cite[Prop.~3.1.2]{ryan_thesis} for another account whose notation matches ours.  The idea is just to do a family version of the argument above. \bigskip

Yoshioka shows that the correspondence $X$ is smooth in \cite[Lem.~5.2]{yoshioka_reflection}, and that the map
\[ f \times g\colon X \to M_{\chi-k} \times M_{\chi+k} \]
is a closed embedding in [ibid., Rmk.~5.1].  The following is related to [ibid., Prop.~2.8]:
\begin{prop}
\label{WZbundles}
Consider the diagram
\[ \begin{tikzcd}
& W \ar[swap]{ld}{\tilde g} \ar{rd}{\tilde f} \ar[dashed]{d}{\pi} & \\
X\ar{d}{e} \ar{dr}[pos=0.3]{f} & Z\ar{dl} \ar{dr} & Y\ar{dl}[swap,pos=0.3]{g} \ar{d}{h} \\
M_{\chi-2m} & M_{\chi} & M_{\chi+2n}
\end{tikzcd} \]
where $X$, $Y$, and $Z$ are the correspondences just constructed and $W$ is the fiber product $X \times_{M_\chi} Y$.
\begin{enumerate}
\renewcommand \labelenumi {(\alph{enumi})}
\item There is a map $\pi\colon W \to Z$ that makes the diagram commute.  It is a $\Gr(n,\,m+n) = \Gr(m+n,\,m)$ bundle.  In particular, $W$ is smooth.
\item The map $\tilde g$ embeds fibers of $\pi$ into fibers of $e$ as sub-Grassmannians, and similarly $\tilde f$ embeds fibers of $\pi$ into fibers of $h$.  Thus $f \circ \tilde g = g \circ \tilde f$ is an embedding on the fibers of $\pi$.
\end{enumerate}
\end{prop}
\begin{proof} A point of $W$ represents a sheaf $F \in M_\chi$, a subspace $\C^m \subset H^0(F)$, and a quotient $H^1(F) \twoheadrightarrow \C^n$.  As we saw earlier, the quotient gives an extension
\[ 0 \to \O_S^n \to G \to F \to 0 \]
for some sheaf $G \in M_{\chi+2n}$, as well as a subspace $\C^n \subset H^0(G)$.  Moreover, the subspace $\C^m \subset H^0(F)$ gives a subspace $\C^{m+n} \subset H^0(G)$ which contains this $\C^n$.  So points in $W$ can be described as flags of subspaces $\C^n \subset \C^{m+n} \subset H^0(G)$, or similarly, as flags of quotients $H^1(E) \twoheadrightarrow \C^{m+n} \twoheadrightarrow \C^m$, where $E \in M_{\chi-2m}$ is defined by the exact sequence
\[ 0 \to \O_S^m \to F \to E \to 0. \]
The map $\pi$ remembers the $\C^{m+n}$ and forgets the $\C^n$ or $\C^m$, so it is a Grassmannian bundle as claimed.  The fiber of $e$ is $\Gr(H^1(E),m)$, and the fiber of $h$ is $\Gr(n,H^0(G))$, so we see that $\tilde g$ and $\tilde f$ embed the fiber of $\pi$ as a sub-Grassmannian.  Now the last claim follows because $e \times f$ and $g \times h$ are closed embeddings.
\end{proof}

\begin{prop} \label{can_use_app}
The results of Appendix \ref{appendix} apply to $f\colon X \to M_{\chi-k}$ if $\chi \le 0$, and to $g\colon X \to M_{\chi+k}$ if $\chi \ge 0$.
\end{prop}
\begin{proof}
First we argue that $R^1 q_* U_{\chi-k}$ and $\sheafExt^2_{q'}(U_{\chi+k},\O)$ have projective dimension 1.  From \cite[Eq.~(70)]{markman_brill_noether} we get an exact sequence
\begin{equation} \label{eq_70}
0 \to R^0 q_* U_{\chi-k} \to V_0 \xrightarrow\rho V_1 \to R^1 q_* U_{\chi-k} \to 0
\end{equation}
on $M_{\chi-k}$, where $V_0$ and $V_1$ are vector bundles twisted by the same Brauer class as $U_{\chi-k}$.  On the one hand, if $\chi \le 0$ then $\chi-k \le 0$, so $R^0 q_* U_{\chi-k}$ is supported on the proper subscheme $_1 M_{\chi-k}$ of $M_{\chi-k}$.  But the vector bundle $V_0$ is torsion-free, so $R^0 q_* U_{\chi-k}$ must vanish, so \eqref{eq_70} gives a two-term resolution of $R^1 q_* U_{\chi-k}$ by vector bundles.  On the other hand, if $\chi \ge 0$ then we write the analogous exact sequence on $M_{\chi+k}$, and dualize to get
\[ 0 \to \sheafExt^1_{q'}(U_{\chi+k},\O) \to V_1'^\vee \to V_0'^\vee \to \sheafExt^2_{q'}(U_{\chi+k},\O) \to 0. \]
By a similar argument we find that $\sheafExt^1_{q'}(U_{\chi+k},\O)$ vanishes, giving a two-term resolution of $\sheafExt^2_{q'}(U_{\chi+k},\O)$ by vector bundles.

Now if $\chi \le 0$ then the Brill--Noether locus $_t M_{\chi-k}$ is the locus where rank of the map $\rho\colon V_0 \to V_1$ drops by $t$, and its codimension is the expected one, which is more than good enough for Appendix \ref{appendix}.  Similarly, if $\chi \ge 0$ then $_t M_{\chi+k}$ is the locus where $V_1'^\vee \to V_0'^\vee$ drops rank by $t$, and it has the expected codimension.
\end{proof}

\subsection*{Allowing negative ranks}
At the beginning of the section we demanded $x > h^2/2$.  Now we allow $x$ to be any integer, and consider the sequence of moduli spaces
\begin{align*}
\dots && M(-1,h,-1-x) && M(0,h,-x) && M(1,h,1-x) && \dots,
\end{align*}
where in the leftward direction we define
\[ M(-r,h,-r-x) := M(r,h,r+x),\footnotemark \]
\footnotetext{As this definition suggests, $M(0,h,-x)$ and $M(0,h,x)$ are isomorphic: both pa\-ra\-me\-trize line bundles or rank-1 torsion-free sheaves on curves in the linear system $|h|$, and the two families of sheaves are related by taking duals.}%
which we refer to hereafter as our ``negative rank fix.''  With this convention, the expected dimension of $M(r,h,r-x)$ is $h^2 - 2r^2 + 2rx + 2$, which is negative for $r \gg 0$ or $r \ll 0$, so only finitely many of the moduli spaces in question are non-empty.

Now we must construct correspondences
\[ \begin{tikzcd}
& X \ar{dl} \ar{dr} & \\ M_{\chi-k} & & M_{\chi+k}.
\end{tikzcd} \]
in three cases.

If both $M_{\chi-k}$ and $M_{\chi+k}$ lie to the right of our negative rank fix, that is, if $M_{\chi-k} = M(r,h,r-x)$ with $r \ge 0$ and thus $M_{\chi+k} = M(r+k,h,r+k-x)$, then we proceed as before: $X$ is the space of pairs consisting of a sheaf $F \in M_{\chi-k}$ and quotient $H^1(F) \twoheadrightarrow \C^k$, or equivalently, a sheaf $G \in M_{\chi+k}$ and a subspace $\C^k \subset H^0(G)$.

If both lie to the left of our negative rank fix, that is, if $M_{\chi+k} = M(-r,h,-r-x)$ with $-r \le 0$ and thus $M_{\chi-k} = M(-r-k,h,-r-k-x)$, then we do the reverse: we let $X$ be the space of pairs consisting of a sheaf $F \in M_{\chi-k}$ and a subspace $\C^k \subset H^0(F)$, or equivalently, a sheaf $G \in M_{\chi+k}$ and a quotient $H^1(G) \twoheadrightarrow \C^k$.

If $M_{\chi-k}$ and $M_{\chi+k}$ straddle our negative rank fix, that is, if $M_{\chi+k} = M(r,h,r-x)$ with $0 < r < k$ and thus $M_{\chi-k} = M(r-k,h,r-k-x)$ with $r-k < 0$, then we proceed as follows.    Again we let $X$ be the space of pairs consisting of a sheaf $G \in M_{\chi+k}$ and a subspace $\C^k \subset H^0(G)$.  But now the evaluation map
\[ \O_S^k \to G \]
cannot be injective, and may not be surjective, so the cone will either be a shifted vector bundle, or a complex with $\H^{-1}$ a vector bundle and $\H^0$ a sheaf with zero-dimensional support.  But if we take the (derived) dual of that cone and shift one place to the left, that is,
\[ F := (\cone(\O_S^k \to G))^\vee[1], \]
then the result is a stable sheaf in $M_{\chi-k}$, as Markman explains in \cite[\S5.8]{markman_brill_noether}.  Going the other way, we can take a sheaf $F \in M_{\chi-k}$ and a subspace $\C^k \subset H^0(F)$, and set $G = (\cone(\O_S^k \to F))^\vee[1]$.  More formally, we let
\[ X := \Gr(\sheafExt^2_q(U_{\chi-k},\O),\,k) \]
with its natural map to $M_{\chi-k}$, or
\[ X := \Gr(\sheafExt^2_{q'}(U_{\chi+k},\O),\,k) \]
with its natural map to $M_{\chi+k}$.  That these two Grassmannians are isomorphic is proved in \cite[Thm.~39]{markman_brill_noether}.

Propositions \ref{WZbundles} and \ref{can_use_app} continue to hold in this more general setting.


\section{An example}

We give an extended example which exhibits the geometric richness typical of the subject.  This section is not logically necessary for the rest of the paper, and we do not give complete proofs of every assertion, but it should clarify some of the subtleties of the construction in the previous section. \bigskip

Let $S$ be a K3 surface of degree 10 and Picard rank 1: thus $S$ is obtained from the Grassmannian $\Gr(2,5) \subset \P^9$ by intersecting with three hyperplanes and a quadric.  The ample generator $h \in \Pic(S)$ embeds $S$ into $\P^6$, which is just the intersection of the three hyperplanes in $\P^9$.

We consider the sequence of moduli spaces
\[ \makebox[\displaywidth]{$\begin{array}{cccccc}
M(-3,h,-2) & M(-2,h,-1) & M(-1,h,0) & M(0,h,1) & M(1,h,2) & M(2,h,3) \\
(\dim = 0) & (\dim = 8) & (\dim = 12) & (\dim = 12) & (\dim = 8) & (\dim = 0)
\end{array}$} \]
We will analyze the correspondences between $M(1,h,2)$ and the other moduli spaces.  We begin by observing that $M(1,h,2)$ is isomorphic to $\Hilb^4(S)$ via the map that sends a length-4 subscheme $\zeta \subset S$ to the twisted ideal sheaf $I_\zeta(h)$.

\subsection*{One step to the left: $M(0,h,1)$}
For a length-4 subscheme $\zeta \subset S$, sections of $I_\zeta(h)$ (up to rescaling) correspond to hyperplanes in $\P^6$ that contain $\zeta$.  Thus in the correspondence
\[ \begin{tikzcd}
& X \ar[swap]{dl}{f} \ar{dr}{g} & \\ M(0,h,1) & & M(1,h,2),
\end{tikzcd} \]
the fiber of $g$ over $I_\zeta(h) \in M(1,h,2)$ is the space of hyperplanes containing $\zeta$.  Four points typically span a $\P^3 \subset \P^6$, so $g$ is generically a $\P^2$ bundle.

A sheaf in $M(0,h,1)$ is typically a line bundle of degree 6 on a curve $C$ in the linear system $|h|$, that is, a curve of the form $S \cap \text{hyperplane}$.  (If $C$ is singular then line bundles may degenerate to rank-1 torsion-free sheaves.)  Such a curve has arithmetic genus 6, and is reduced and irreducible thanks to our assumption that $\Pic(S) = \Z h$.  Thus $M(0,h,1)$ is an example of a Beauville--Mukai integrable system: it is fibered in Abelian 6-folds over $|h| \cong \P^6$, and more precisely is the compactified relative $\Pic^6$ of the universal family of curves over that linear system.

To study the map $f$, take a point of $X$, that is, a length-4 subscheme $\zeta \subset S$ and a hyperplane $H$ that contains it.  We get a section
\[ \O_S \to I_\zeta(h), \]
and we find that the cokernel is $\omega_C(-\zeta) \in M(0,h,1)$, where again $C = S \cap H$.  But $f$ is not surjective: by Serre duality we have $h^1(\omega_C(-\zeta)) = h^0(\O_C(\zeta))$ which is positive because $\zeta$ is an effective divisor, whereas a general degree-6 line bundle on $C$ has $h^1 = 0$.  In fact $f$ is generically injective, mapping $X$ birationally onto the first Brill--Noether stratum of $M(0,h,1)$, which we see has codimension 2: we have $\dim M(1,h,2) = 8$, so $\dim X = 10$, but $\dim M(0,h,1) = 12$.

\subsection*{One step to the right: $M(2,h,3)$}
Four points $\zeta \subset S$ typically span a $\P^3$, but they might only span a plane, giving $h^0(I_\zeta(h)) = 4$ and $h^1(I_\zeta(h)) = 1$.  Using Serre duality $H^1(I_\zeta(h)) \cong \Ext^1(I_\zeta(h), \O_S)^*$ we get an extension
\[ 0 \to \O_S \to E \to I_\zeta(h) \to 0, \]
where $E \in M(2,h,3)$.  The latter moduli space is a single point, so it is not hard to guess that $E = T^\vee|_S$, where $T \subset \O_{\Gr}^5$ is the tautological sub-bundle on $\Gr(2,5)$.

Going the other way, from $M(2,h,3)$ back to $M(1,h,2)$, we find that $h^0(T^\vee|_S) = h^0(T^\vee) = 5$, so our correspondence is a $\P^4$:
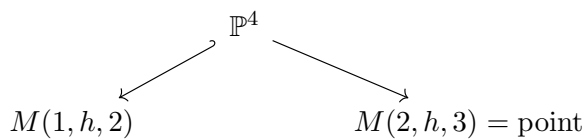
\[ \begin{tikzcd}
& \P^4 \ar[dl,hook'] \ar[dr] & \\ M(1,h,2) & & M(2,h,3) = \text{point}
\end{tikzcd} \]
We find that $c_2(T^\vee|_S) = 4$, so a non-zero section of $T^\vee|_S$ vanishes at 4 points $\zeta \subset S$, and the cokernel of the section $\O_S \to T^\vee|_S$ must be $I_\zeta(h)$.  The map from the correspondence to $M(1,h,2)$ is injective, giving a Lagrangian $\P^4$ in a holomorphic symplectic 8-fold.

\subsection*{Can't go two steps to the right: $M(3,h,4)$}
We might think that four points in $S$ could degenerate even further and only span a line, but this is impossible, as we can see in two ways.

For a classical argument, we know that $\Gr(2,5)$ is an intersection of quadrics in $\P^9$, so the K3 surface $S$ is an intersection of quadrics in $\P^6$.  If any line meets $S$ in four points then it meets each quadric in four points, hence is contained in each quadric, hence is contained in $S$, but this contradicts our assumption that $\Pic(S) = \Z h$.

For an argument using the present moduli set-up, we can say that if $h^0(I_\zeta(h)) \ge 5$ then $h^1(I_\zeta(h)) \ge 2$, giving an extension of $I_\zeta(h)$ by $\O_S^2$ and hence a point of $M(3,h,4)$.  But the expected dimension of this moduli space is $-12$, so it is empty.

\subsection*{Two steps to the left: $M(-1,h,0) := M(1,h,0)$}
A 2-dimensional subspace of $H^0(I_\zeta(h))$ corresponds to a $\P^4 \subset \P^6$ that contains $\zeta$.  Thus in the correspondence
\[ \begin{tikzcd}
& Y \ar[swap]{dl}{e} \ar{dr}{k} & \\ M(-1,h,0) & & M(1,h,2),
\end{tikzcd} \]
the map $k$ is generically a $\P^2$ bundle: a general $\zeta \in \Hilb^4(S)$ spans a $\P^3$, and the space of $\P^4$s that fit between this $\P^3$ and the ambient $\P^6$ is a $\P^2$.  This fiber is dual to the $\P^2$ fiber of $g\colon X \to M(1,h,2)$ seen earlier, which was the space of hyperplanes that fit between $\P^3$ and $\P^6$.  Over the Brill--Noether locus of $M(1,h,2)$ where $\zeta$ only spans a plane, $k$ becomes a $\Gr(2,4)$ bundle, whereas $g$ became a $\P^3$ bundle.

Our negative rank fix defines $M(-1,h,0)$ to mean $M(1,h,0)$, which is isomorphic to $\Hilb^6(S)$ via the map that sends a length-6 subscheme $\eta$ to the twisted ideal sheaf $I_\eta(h)$.

To study the map $e$, take a point of $Y$, that is, a length-4 subscheme $\zeta \subset S$ and a $\P^4$ that contains it.  Then $\xi := S \cap \P^4$ must be a length-10 subscheme: it cannot be a curve by our assumption that $\Pic(S) = \Z h$, so it must be a finite subscheme, and the embedding $S \subset \P^6$ has degree 10.  We would like to say that $e$ sends our point of $Y$ to $I_\eta(h) \in M(1,h,0)$, where $\eta = \xi \setminus \zeta$.  This is valid if $\xi$ is reduced, but in general the difference $\xi \setminus \zeta$ is not well-defined, so we must be a little more careful.  Take the surjection $\O_\xi \twoheadrightarrow \O_\zeta$; dualize to get an injection $\omega_\zeta \hookrightarrow \omega_\xi$; observe that $\omega_\xi = \O_\xi$, because $\xi$ is a complete intersection in $S$; thus $\omega_\xi/\omega_\zeta$ is a quotient of $\O_\xi$, so it is $\O_\eta$ for some subscheme $\eta \subset \xi$, whose length we find to be 6.  Markman discusses this construction in \cite[Example 40]{markman_brill_noether}.

Our recipe from the end of \S\ref{sectioncorrespondences} said to take a 2-dimensional subspace of $H^0(I_\zeta(h))$, take the cone on the associated map $\O_S^2 \to I_\zeta(h)$, take the derived dual, and shift by 1.  The energetic reader may check that this really produces $I_\eta(h)$, using the Koszul resolution
\[ 0 \to \O_S(-2h) \to \O_S(-h)^2 \to I_\xi \to 0 \]
and the exact sequence
\[ 0 \to I_\xi \to I_\zeta \to \omega_\eta \to 0. \]

The map $e$ is generically injective, but not surjective: its image consists of length-6 subschemes contained in a $\P^4$, whereas six points would typically span a $\P^5 \subset \P^6$.

\subsection*{Three steps to the left: $M(-2,h,-1) := M(2,h,1)$}
We have seen that for any $\zeta \in \Hilb^4(S)$, either $h^0(I_\zeta(h)) = 3$, which is the generic behavior, or $h^0 = 4$, which occurs along a Lagrangian $\P^4$.  In the first case, we take the cone of
\[ \O_S^3 \to I_\zeta(h), \]
take the derived dual, and shift by 1 to get a rank-2 stable sheaf in $M(2,h,1)$.  We do not know any more down-to-earth description of the latter moduli space, but we remark that the operation just described is a birational map, in fact the Mukai flop of $M(1,h,2)$ along the Lagrangian $\P^4$.  Our correspondence between $M(2,h,1)$ and $M(1,h,2)$ is the graph of that birational map, blowing up the $\P^4$ and blowing down the exceptional divisor in the other direction.

\subsection*{Four steps to the left: $M(-3,h,-2) := M(3,h,2)$}
Suppose that $\zeta$ is in the Brill--Noether locus of $\Hilb^4(S)$; that is, it only spans a plane in $\P^6$, so $h^0(I_\zeta(h)) = 4$.  We have seen that $I_\zeta(h)$ is a quotient of $T^\vee|_S$, so it is globally generated.  Thus the kernel of
\[ \O_S^4 \twoheadrightarrow I_\zeta(h) \]
is a stable vector bundle; in fact it is $Q^\vee|_S$, where $Q = \O_{\Gr}^5/T$ is the tautological rank-3 quotient bundle on the Grassmannian.  Thus the cone of the surjection above is $Q^\vee|_S[1]$, so its dual is $Q|_S[-1]$, and shifting by 1 we get $Q|_S$ which is the unique point of $M(3,h,2)$, as expected.

\subsection*{Other remarks}
The space $M(1,h,0) \cong \Hilb^6(S)$ has three Brill-Noether strata, depending on whether the 6 points span a hyperplane in $\P^6$, or a $\P^4$, or a $\P^3$.  They cannot span a plane, for then every quadric containing $S$ would intersect that plane in the same conic (no four of the points can be collinear, as we saw earlier, so through any five of them there passes a unique conic) and thus $S$ would contain a conic, contradicting our hypothesis that $\Pic(S) = \Z h$.

The Abelian-fibered $M(0,h,1)$ also has three Brill--Noether strata.  The correspondence between $M(1,h,0)$ and $M(0,h,1)$ is the graph of a birational map described by Beauville in \cite[Prop.~1.3]{beauville_counting}.  It sends a length-6 subscheme $\xi \subset S$ that spans a hyperplane $H \subset \P^6$ to $\O_C(\xi)$, where $C$ is the curve $S \cap H$.  The analogous birational map for lower-degree K3 surfaces was studied in terms of derived categories in \cite{adm}.

We have seen that $M(2,h,3)$ is a single point representing $T^\vee|_S$, and $M(3,h,2)$ is a single point representing $Q|_S$.  The 5-step correspondence between them amounts to the tautological exact sequence
\[ 0 \to T|_S \to \O_S^5 \to Q|_S \to 0. \]


\section{Line bundles}
\label{sectionlinebundles}

In \S\ref{sectioncorrespondences} we constructed correspondences
\begin{equation} \label{kstep}
\begin{tikzcd}
& X \ar[swap]{dl}{f} \ar{dr}{g} & \\ M_{\chi-k} & & M_{\chi+k}
\end{tikzcd}
\end{equation}
for each $\chi$ and $k$ such that $X_{\chi \pm k}$ are not empty.  Now we will choose line bundles $L_X$ on each $X$, and in Definition \ref{maindef} set
\begin{equation} \label{EF_def}
\begin{array}{rll}
E^{(k)}_\chi &= (f \times g)_* L_X & \in D(M_{\chi-k} \times M_{\chi+k}), \text{ and} \\
F^{(k)}_\chi &= (g \times f)_* (L_X^{-1} \otimes \omega_X) & \in D(M_{\chi+k} \times M_{\chi-k}).
\end{array}
\end{equation}

In fact we have quite a bit of freedom in our choice of line bundles:\footnote{Sabin Cautis pointed out to us that this degree of freedom may represent an interesting enrichment of the $\sl_2$ structure, just as the $\C^\times$ action in his example actually allows for a categorification of the quantum group $U_q(\sl_2)$ \cite[Rmk.~3.1]{ckl1795}.}
the only requirement is that for every diagram as in Proposition \ref{WZbundles},
\[ \begin{tikzcd}
& W \ar[swap]{ld}{\tilde g} \ar{rd}{\tilde f} \ar{d}{\pi} & \\
X\ar{d}{e} \ar{dr}[pos=0.3]{f} & Z\ar{dl} \ar{dr} & Y\ar{dl}[swap,pos=0.3]{g} \ar{d}{h} \\
M_{\chi-2m} & M_{\chi} & M_{\chi+2n}
\end{tikzcd} \]
we should have
\begin{equation} \label{requirement}
\tilde g^*L_X \otimes \tilde f^*L_Y = \pi^*L_Z \otimes \omega_\pi.
\end{equation}

It follows from Proposition \ref{excess_prop} below that
\begin{equation} \label{excess_omega}
\tilde g^* \omega_X \otimes \tilde f^* \omega_Y = \pi^* \omega_Z \otimes \omega_\pi^2,
\end{equation}
so our first idea is that $L_X$ might be a square root of $\omega_X$.  But this does not work in general: in the previous section's example, we saw that the correspondence ``one step to the right'' is a $\P^4$ bundle over $M(1,h,2)$, so $\omega_X$ has degree $-5$ on the fibers, and in particular $\omega_X$ has no square root.  With a little adjustment, however, the idea can be made to work:

\begin{thm} \label{lb_thm}
For each $\chi$ there is a line bundle $L_\chi \in \Pic(M_\chi)$, such that on any correspondence as in \eqref{kstep}, the line bundle
\[ \omega_X \otimes f^* L_{\chi-k} \otimes g^* L_{\chi+k}^{-1} \in \Pic(X) \]
has a square root, which we call $L_X$.\footnote{We are sorry that $L_X$ and $L_\chi$ look so similar.}
\end{thm}
The Picard groups of our moduli spaces are torsion-free, and the same is true of our correspondences by Proposition \ref{app:ik+2}, so if the promised square root exists then it is unique.  Moreover, \eqref{excess_omega} implies \eqref{requirement}: pulling everything back to $W$ and suppressing pullbacks from the notation, we have
\begin{align*}
L_X^2 \otimes L_Y^2
&= (\omega_X \otimes L_{\chi-2m} \otimes L_\chi^{-1}) \otimes (\omega_Y \otimes L_\chi \otimes L_{\chi+2n}^{-1}) \\
&= (\omega_X \otimes \omega_Y) \otimes (L_{\chi-2m} \otimes L_{\chi+2n}^{-1}) \\
&= (\omega_Z \otimes \omega_\pi^2) \otimes (L_{\chi-2m} \otimes L_{\chi+2n}^{-1}) \\
&= L_Z^2 \otimes \omega_\pi^2.
\end{align*}
Because $W$ is a Grassmannian bundle over $Z$, its Picard group is torsion-free as well, so this implies \eqref{requirement}. \bigskip

The proof of Theorem \ref{lb_thm} will occupy the rest of this section.  We begin by reviewing the standard description of the Picard groups of our moduli spaces, due to O'Grady \cite{og}.

In \S\ref{sectionbn} we reviewed the the Mukai vector of a sheaf or complex on the K3 surface $S$; now we recall the Mukai pairing,
\[ \langle  (a,bh,c),\ (a',b'h,c') \rangle = bb' h^2 - ac' - a'c, \]
which is cooked up so that for two sheaves $F$ and $F'$ we have
\[ \chi(F,F') := \sum (-1)^i \dim \Ext^i_S(F,F') = -\langle v(F), v(F') \rangle, \]
ultimately by Riemann--Roch.

For a moduli space $M = M(r,h,s)$, the \emph{Mukai map} $\theta$ from
\[ (r,h,s)^\perp \subset H^*(S,\Z) \]
to $\Pic(M)$ is defined as follows.  Let $S \xleftarrow{p} S \times M \xrightarrow{q} M$ be the two projections, and choose a (possibly twisted) universal sheaf $U$ on $S \times M$.  For a vector $(a,bh,c) \in (r,h,s)^\perp$, choose a complex $E \in D(S)$ with $v(E) = (a,bh,c)$, and define
\[ \theta(a,bh,s) = \det q_*(U \otimes p^* E^\vee). \]
The complex $q_*(U \otimes p^* E^\vee)$ has rank zero, so if $U$ was twisted by $q^* \alpha$ for some Brauer class $\alpha \in \Br(M)$ then the determinant is twisted by $\alpha^0$, that is, it is naturally untwisted.  Relatedly, if we choose a different universal sheaf $U' = U \otimes q^* L$ for some $L \in \Pic(M)$, then the determinant is tensored by $L^0$, that is, it doesn't change.  In fact the determinant line bundle does not depend on the choice of complex $E$, only on its Mukai vector $(a,bh,c)$, and O'Grady showed that the Mukai map is surjective, and indeed injective if $\dim M > 2$. \bigskip

Writing $M_\chi = M(r,h,s)$, and letting $\theta_\chi$ be the Mukai map for $M_\chi$, the line bundle $L_\chi$ that we will choose is
\[ L_\chi := \theta_\chi(r,0,-s). \]
Observe that
\[ \langle (r,0,-s),\,(r,h,s) \rangle = 0 - rs + rs = 0, \]
so $(r,0,-s)$ is in the domain of the Mukai map. \bigskip

To prove Theorem \ref{lb_thm}, we need to extend the domain of the Mukai map, allowing vectors $(a,bh,c)$ for which the pairing $\langle (a,bh,c),\,(r,h,s) \rangle =: n$ is not zero.  Then the complex $q_*(U \otimes p^* E^\vee)$ has rank $-n$, so if the universal sheaf $U$ is twisted by $q^* \alpha$ then the determinant line bundle is twisted by $\alpha^{-n}$.  The reader may object that a twisted line bundle can always be untwisted, and we admit that $\alpha^{-n}$ is trivial in $\Br(M)$, but it is not canonically trivial: if we represent $\alpha$ by a \v Cech 2-cocycle, then $\alpha^{-n}$ is the coboundary of some 1-cochain, but different choices of 1-cochain will give untwisted sheaves that differ by a line bundle.  This would ruin our calculations, so we will take care to avoid unnatural untwisting.  Relatedly, if $n \ne 0$ then $\theta(a,bh,c)$ depends on our choice of universal bundle: if we choose another one $U' = U \otimes q^* L$ for some $L \in \Pic(M)$, then $\theta(a,bh,c)$ changes by $L^{-n}$. \bigskip

Now consider the correspondence
\[ \begin{tikzcd}
& X \ar[swap]{dl}{f} \ar{dr}{g} & \\ M_{\chi-k} & & M_{\chi+k}.
\end{tikzcd} \]
Assume for simplicity that $\chi \ge 0$ and that both moduli spaces lie to the right of the negative rank fix; the other cases are similar.  Fix universal sheaves $U_{\chi\pm k}$ on $S \times M_{\chi\pm k}$, twisted by Brauer classes $\alpha_{\chi\pm k} \in \Br(M_{\chi\pm k})$.

Because $\chi \ge 0$, the map $g$ is surjective, and $f$ is generically injective.  In \S\ref{sectioncorrespondences} we defined
\[ X = \Gr(\sheafExt^2_{q'}(U_{\chi+k},\O_S),\,k). \]
We have a universal quotient
\begin{equation} \label{univ_quot}
g^* \sheafExt^2_{q'}(U_{\chi+k},\O_S) \twoheadrightarrow Q,
\end{equation}
where $Q$ is a rank-$k$ vector bundle on $X$, twisted by $g^* \alpha_{\chi+k}^{-1}$.  By Proposition \ref{app:ik}(b) we have
\[ \omega_X = g^*(\det \sheafExt^2_{q'}(U_{\chi+k},\O_S))^k \otimes (\det Q)^{-\chi-k}. \]
We can rewrite this as
\begin{align}
\omega_X &= g^*(\det q'_* U_{\chi+k})^{-k} \otimes (\det Q)^{-\chi-k} \notag \\
&= g^* \theta_{\chi+k}(1,0,1)^{-k} \otimes (\det Q)^{-\chi-k} \notag \\
&= g^* \theta_{\chi+k}(-k,0,-k) \otimes (\det Q)^{-\chi-k}, \label{omega_X}
\end{align}
where in the first line we have used Grothendieck duality, and in the second and third lines we have embraced our extended definition of the Mukai map.  We are reassured to see that $\omega_X$ is naturally untwisted: the first factor of \eqref{omega_X} is twisted by $g^* \alpha_{\chi+k}^{-k(\chi+k)}$, and the second is twisted by $g^* \alpha_{\chi+k}^{-k (-\chi-k)}$, so the twists cancel. \bigskip

Now we want to compare the line bundles $L_{\chi\pm k}$ when pulled back to $X$, so we need to compare the universal sheaves $U_{\chi\pm k}$ when pulled back to $S \times X$.\footnote{See \cite[\S2.4]{yoshioka_reflection} for a related calculation.}  The quotient \eqref{univ_quot} gives rise to an injection
\[ \O_S \boxtimes Q^\vee \hookrightarrow (1 \times g)^* U_{\chi+k} \]
whose cokernel is an $X$-flat family of sheaves on $S$ that determines the map $f\colon X \to M_{\chi-k}$.  Thus there is a line bundle $L$ on $X$ and an exact sequence
\[ 0 \to \O_S \boxtimes Q^\vee \to (1 \times g)^* U_{\chi+k} \to (1 \times f)^* U_{\chi-k} \otimes \pi_X^* L \to 0. \]
Write $M_{\chi-k} = M(r,h,s)$, so $L_{\chi-k}$ is $\theta_{\chi-k}(r,0,-s)$.  Choose an $E \in D(S)$ with $v(E) = (r,0,-s)$.  Tensor the exact sequence above with $E^\vee$ and push down to $X$.  The first term becomes $\RGamma(E^\vee) \otimes Q^\vee$, whose determinant is $(\det Q)^{-\chi(E^\vee)} = (\det Q)^{-r+s}$.  The second term becomes a complex whose determinant is $g^* \theta_{\chi+k}(r,0,-s)$.  The third term becomes a complex whose determinant is $f^* \theta_{\chi-k}(r,0,-s) \otimes L^0$.  Thus as twisted line bundles on $X$ we have
\[ (\det Q)^{-r+s} \otimes f^* \theta_{\chi-k}(r,0,-s) = g^* \theta_{\chi+k}(r,0,-s). \]
We also have $M_{\chi+k} = M(r+k,h,s+k)$, so
\begin{align}
f^* L_{\chi-k} \otimes g^* L_{\chi+k}^{-1}
&= f^* \theta_{\chi-k}(r,0,-s) \otimes g^* \theta_{\chi+k}(r+k,0,-s-k) \notag \\
&= g^* \theta_{\chi+k}(2r+k,0,-2s-k) \otimes (\det Q)^{r-s}. \label{L_minus_L}
\end{align}

Finally we are in a position to take the square root of
\[ \omega_X \otimes f^* L_{\chi-k} \otimes g^* L_{\chi+k}^{-1}. \]
Tensoring \eqref{omega_X} and \eqref{L_minus_L}, we see that this equals
\[ g^* \theta_{\chi+k}(2r,0,-2s-2k) \otimes (\det Q)^{-2s-2k},\]
so the desired square root is
\[ L_X := g^* \theta_{\chi+k}(r,0,-s-k) \otimes (\det Q)^{-s-k}. \]
We check that this is naturally untwisted: the first factor is twisted by $g^* \alpha_{\chi+k}^{-k(s+k)}$, and the second is twisted by $g^* \alpha_{\chi+k}^{-k(-s-k)}$.


\section{Straightforward conditions}
\label{sectioneasycond}

Now we turn to the task of verifying that our construction satisfies the conditions of Definition \ref{maindef}.  In this section we check the easy conditions (i), (ii), (vi), and (vii).  In subsequent sections we check the harder conditions (iii) and (v), and finally the hardest condition (iv).

\renewcommand \thecond {\ref{condi}}
\begin{cond}
Each Hom space between two objects of $D(M_\chi)$ is finite dimensional.
\end{cond}
\begin{proof}
This is because $M_\chi$ is smooth and proper.
\end{proof}

\renewcommand \thecond {\ref{condii}}
\begin{cond}
The right adjoint of $E^{(k)}_\chi$ is $F^{(k)}_\chi[k\chi]$, and the left adjoint is $F^{(k)}_\chi[-k\chi]$.
\end{cond}
\begin{proof}
Consider the correspondence
\[ \begin{tikzcd}
& X \ar[swap]{dl}{f} \ar{dr}{g} & \\ M_{\chi-k} & & M_{\chi+k}
\end{tikzcd} \]
constructed in \S\ref{sectioncorrespondences}, and the line bundle $L_X \in \Pic(X)$ constructed in Theorem \ref{lb_thm}.  From our definitions of $E^{(k)}_\chi$ and $F^{(k)}_\chi$ as kernels in \eqref{EF_def}, we see that as functors,
\[ E^{(k)}_\chi = g_* (L_X \otimes f^*(-))\colon D(M_{\chi-k}) \to D(M_{\chi+k}), \]
and
\[ F^{(k)}_\chi = f_* (L_X^{-1} \otimes \omega_X \otimes g^*(-))\colon D(M_{\chi+k}) \to D(M_{\chi-k}), \]
using \cite[Ex.~5.12]{huybrechts_fm}.  Thus the right adjoint of $E^{(k)}_\chi$ is
\begin{align*}
f_* (L_X^{-1} \otimes g^!(-)) 
&= f_* (L_X^{-1} \otimes g^*(-) \otimes \omega_X[\dim g]) \\
&= F^{(k)}_\chi[\dim g],
\end{align*}
where in the first equality we have used the fact that the canonical bundle of $M_{\chi+k}$ is trivial, and the left adjoint is
\begin{align*}
f_! (L_X^{-1} \otimes g^*(-))
&= f_* (L_X^{-1} \otimes g^*(-) \otimes \omega_X[\dim f]) \\
&= F^{(k)}_\chi[\dim f].
\end{align*}
It remains to show that $\dim g = k\chi$ and $\dim f = -k\chi$.

If $\chi \geq 0$ then $g$ is generically a $\Gr(k,\chi+k)$ bundle, so $\dim g = k\chi$ as desired.  For $\dim f$, it is quickest to compute the dimensions of the moduli spaces: if $M_{\chi+k} = M(r,h,s)$ with $r+s = \chi+k$, then $M_{\chi-k} = M(r-k,h,s-k)$, and we have
\begin{align*}
\dim f &= \dim X - \dim M_{\chi-k} \\
&= \dim M_{\chi+k} + k\chi - \dim M_{\chi-k} \\
&= (h^2 - 2rs + 2) + k\chi - (h^2 - 2(r-k)(s-k) + 2)
\end{align*}
which simplifies to $-k\chi$ as desired.

If $\chi \le 0$ then $f$ is generically a $\Gr(k,|\chi|+k)$ bundle, and the calculation is similar.
\end{proof}

\renewcommand \thecond {\ref{condvi}}
\begin{cond}
If $\chi \leq 0$ then for $k' > k$, the image of $\supp(E^{(k)}_{\chi-k})$ under the projection to $M_\chi$ is not contained in the image of $\supp(E^{(k')}_{\chi-k'})$ also under the projection to $M_\chi.$ Similarly, if $\chi \geq 0$ then for $k' > k$, the image of $\supp(E^{(k)}_{\chi+k})$ in $M_\chi$ is not contained in the image of $\supp(E^{(k')}_{\chi+k'})$.
\end{cond}
\begin{proof}
In either case, the image of $\supp(E^{(k)})$ in $M_\chi$ is the Brill--Noether locus $_k M_\chi$, which is strictly bigger than $_{k'} M_\chi$.
\end{proof}

\renewcommand \thecond {\ref{condvii}}
\begin{cond}
All $E^{(r)}$s and $F^{(r)}$s are sheaves.
\end{cond}
\begin{proof}
This is because $f \times g \colon X \to M_{\chi-k} \times M_{\chi+k}$ is a closed embedding.
\end{proof}


\section{Divided power condition (iii)}
\label{sectioncond3}

We will now see that our line bundles $L_X$ in \S\ref{sectionlinebundles} were constructed precisely to give condition (iii) in Definition \ref{maindef}.  In fact it is no harder, and maybe a little clearer, to prove a slightly more general statement:
\begin{thm}[Implies condition \ref{condiii}] \label{thmiii}
At the level of cohomology of complexes we have
\[\H^*(E^{(n)}_{\chi+n} \ast E^{(m)}_{\chi-m}) \cong E^{(m+n)}_{\chi-m+n} \otimes_\C H^*(\Gr(n,m+n)),\]
where the grading of $H^*(\Gr(n,m+n))$ is centered around 0.
\end{thm}

Recall the set-up of \S\ref{sectionlinebundles}: we have a diagram
\[ \begin{tikzcd}
& W \ar[swap]{dl}{\tilde g} \ar{d}{\pi} \ar{dr}{\tilde f} & \\
X \ar[swap]{d}{e} \ar[pos=0.3]{dr}{f} & Z \ar{dl}[swap, pos=0.3]{k} \ar[pos=0.3]{dr}{\ell} & Y \ar[swap, pos=0.3]{dl}{g} \ar{d}{h} \\
M_{\chi-2m} & M_\chi & M_{\chi+2n}
\end{tikzcd} \]
where $W$ is the fiber product of $X$ and $Y$ over $M_\chi$, $\pi$ is a $\Gr(n,m+n)$ bundle, and we have line bundles $L_X$, $L_Y$, and $L_Z$ on the correspondences that satisfy
\[ \tilde g^*L_X \otimes \tilde f^*L_Y = \pi^*L_Z \otimes \omega_\pi. \]

As a functor, the composition $E^{(n)}_{\chi+n} \ast E^{(m)}_{\chi-m}$
is
\[ h_*(L_Y \otimes g^* f_*(L_X \otimes e^*(-))). \]
We cannot expect the cohomology-and-base-change map $g^* f_* \to \tilde f_* \tilde g^*$ to be an isomorphism, because the dimension of $W$ is bigger than expected.\footnote{If $\chi \ge m$, for example, then $g$ is birational onto its image $_n M_\chi \subset M_\chi$, which has codimension $n(\chi+n)$, and $f$ is generically a $\Gr(m,\chi)$ bundle; but over the image of $g$, the map $f$ is generically a $\Gr(m,\chi+n)$ bundle, so the actual dimension of $W$ exceeds the expected dimension $\dim X + \dim Y - \dim M_\chi$ by $\dim \Gr(m,\chi+n) - \dim \Gr(m,\chi) = mn$.  The other two cases $-n \le \chi \le m$ and $\chi \le -n$ are similar.}
Instead we must study the \emph{excess normal bundle} $E$ of this fiber square, defined by the exact sequence
\begin{equation} \label{excess}
0 \to T_W \to \tilde g^* T_X \oplus \tilde f^* T_Y \to \tilde g^* f^* T_{M_\chi} \to E \to 0.
\end{equation}
The bundle $E$ measures the failure of $f$ and $g$ to be transverse in the sense of manifolds.  In Proposition \ref{excess_prop} below, we will see that $E$ is the relative cotangent bundle $\Omega_\pi$, but first let us see how that implies Theorem \ref{thmiii}.

\begin{proof}[Proof of Theorem \ref{thmiii}]
Switching back to the language of Fourier--Mukai kernels, the convolution $E^{(n)}_{\chi+n} \ast E^{(m)}_{\chi-m}$ can be obtained in three steps:

\paragraph{Step 1:}  Start with $\O_\Delta \in D(M_\chi \times M_\chi)$, and apply $(f \times g)^*$ to get a kernel in $D(X \times Y)$ that induces the functor $g^* f_*$. \bigskip

For the fiber square
\[ \begin{tikzcd}
W \ar[hookrightarrow]{r} \ar[swap]{d}{\tilde g \circ f} & X \times Y \ar{d}{f \times g} \\
\Delta \ar[hookrightarrow]{r} & M_\chi \times M_\chi,
\end{tikzcd} \]
we consider the excess normal bundle in the sense of Fulton \cite[\S6.3]{fulton}, defined by the exact sequence
\[ 0 \to N_{W/X \times Y} \to N_{\Delta/M_\chi \times M_\chi} \to E \to 0, \]
which measures the failure of $f \times g$ to be transverse to $\Delta$.  Note that both inclusions are regular embeddings because the spaces involved are smooth.  To see that this excess normal bundle coincides with the one defined in \eqref{excess}, observe that $N_{W/X \times Y}$ is the cokernel of the first map in \eqref{excess}, and $N_{\Delta/M_\chi \times M_\chi} = T_{M_\chi}$. \smallskip 

Now by \cite[Lem.~3.2]{thomason},\footnote{We thank Adeel Khan for this reference.} we have
\[ L_j (f \times g)^* \O_\Delta = \Lambda^j E^* = \Lambda^j T_\pi. \]
We do not claim that the complex $(f \times g)^* \O_\Delta$ is the pushforward of a complex on $W$, only that its cohomology sheaves are pushed forward from $W$; there may be extensions between these as sheaves on $X \times Y$ that do not come from $W$.

\paragraph{Step 2:} Tensor $(f \times g)^* \O_\Delta$ with $L_X \boxtimes L_Y$ to get a kernel $K \in D(X \times Y)$ that induces the functor $L_Y \otimes g^* f_*(L_X \otimes -)$. \bigskip

On the level of cohomology sheaves, this amounts to tensoring with $\tilde g^* L_X \otimes \tilde f^* L_Y$ on $W$.  By construction this equals $\pi^* L_Z \otimes \omega_\pi$, so
\begin{align*}
\H^q(K) &= \pi^* L_Z \otimes \omega_\pi \otimes \Lambda^{-q} T_\pi \\
&= \pi^* L_Z \otimes \Omega^{mn+q}_\pi.
\end{align*}

\paragraph{Step 3:} Apply $(e \times h)_*$ to get $E^{(n)}_{\chi+n} \ast E^{(m)}_{\chi-m}$ in $D(M_{\chi-2m} \times M_{\chi+2n})$. \bigskip

We use the Grothendieck spectral sequence
\[ E_2^{p,q} = R^p (e \times h)_* \H^q(K) \Longrightarrow \H^{p+q}((e \times h)_* K). \]
Because the cohomology sheaves $\H^q(K)$ are supported on $W$, we see that $(e \times h)_*$ acts on them as $(k \times \ell)_* \circ \pi_*$.  We calculate
\begin{align*}
R^p \pi_* \H^q(K) 
&= R^p \pi_*(\pi^* L_Z \otimes \Omega^{mn+q}_\pi) \\
&= L_Z \otimes R^p \pi_*(\Omega^{mn+q}_\pi) \\
&= L_Z \otimes_{\C} H^{mn+q,p}(\Gr(m,m+n)),
\end{align*}
where in the second line we have used the projection formula and in the third we have used the fact that $\pi$ is a $\Gr(m,m+n)$ bundle.  Now $H^{mn+q,p}(\Gr)$ is zero if $mn+q \ne p$, and is $H^{2p}(\Gr)$ if $mn+q=p$.  Thus the spectral sequence degenerates at the $E_2$ page and there are no possible extensions.

We recall that $E^{(m+n)}_{\chi-m+n}$ was defined to be $(k \times \ell)_* L_Z$, so this gives Theorem \ref{thmiii}.
\end{proof}

\begin{prop} \label{excess_prop}
In the exact sequence \eqref{excess}, we have $E \cong \Omega_\pi$.  Moreover, if we use the holomorphic symplectic form $\sigma_\chi$ on $M_\chi$ to identify $T_{M_\chi}$ with $\Omega_{M_\chi}$, then the map $\tilde g^* f^* T_{M_\chi} \to E$ agrees with the restriction map $\tilde g^* f^* \Omega_{M_\chi} \to \Omega_\pi$.
\end{prop}
\begin{proof}
First we observe that the restriction map $\tilde g^* f^* \Omega_{M_\chi} \to \Omega_\pi$ is surjective; this follows from Proposition \ref{WZbundles}(b).

Next we claim that the composition
\[ \tilde g^* T_X \oplus \tilde f^* T_Y \to \tilde g^* f^* T_{M_\chi} \xrightarrow{\sigma_\chi} \tilde g^* f^* \Omega_{M_\chi} \to \Omega_\pi \]
is zero.  It is equivalent to argue that for a point $p \in W$ and tangent vectors $u \in T_{X,\tilde g(p)}$, $v \in T_{Y,\tilde f(p)}$, and $w \in T_{\pi,p}$, we have
\[ \sigma_\chi(f_* u + g_* v,\, f_* \tilde g_* w) = 0. \]
Rewrite the left-hand side as follows:
\begin{align*}
& \sigma_\chi(f_* u,\, f_* \tilde g_* w) + \sigma_\chi(g_* v,\, g_* \tilde f_* w) \\
&= (f^* \sigma_\chi)(u,\, \tilde g_* w) + (g^* \sigma_\chi)(v,\, \tilde f_* w) \\
&= (e^* \sigma_{\chi-2m})(u,\, \tilde g_* w) + (h^* \sigma_{\chi+2n})(v,\, \tilde f_* w) \\
&= \sigma_{\chi-2m}(e_* u,\, e_* \tilde g_* w) + \sigma_{\chi+2n}(h_* v,\, h_* \tilde f_* w) \\
&= \sigma_{\chi-2m}(e_* u,\, k_* \pi_* w) + \sigma_{\chi+2n}(h_* v,\, \ell_* \pi_* w),
\end{align*}
where in the third line we have used Proposition \ref{symp_forms} below.  Now $\pi_* w = 0$ because $w$ is tangent to a fiber of $\pi$, so everything vanishes as desired.

Now the sequence of vector spaces
\[ 0 \to T_{W,w} \to T_{X,\tilde g(w)} \oplus T_{Y,\tilde f(w)} \to T_{M_\chi, f(\tilde g(w))} \to \Omega_{\pi,w} \to 0 \]
is a complex, and the dimensions are right to make it exact, so the corresponding sequence of vector bundles is exact.  The first three terms are the first three terms of \eqref{excess}, so the last terms agree, as desired.
\end{proof}

\begin{prop} \label{symp_forms}
Consider the correspondence
\[ \begin{tikzcd}
& X \ar[swap]{dl}{f} \ar{dr}{g} & \\ M_{\chi-k} & & M_{\chi+k}.
\end{tikzcd} \]
Let $\sigma_{\chi\pm k}$ be the holomorphic symplectic forms on $M_{\chi\pm k}$.  Then $f^* \sigma_{\chi-k}$ and $g^* \sigma_{\chi+k}$ agree up to a non-zero constant.
\end{prop}
\begin{proof}
This is closely related to a result of Mukai on coisotropic $\P^n$ bundles \cite[Prop.~3.1]{mukai_inventiones}.

Suppose that $\chi \ge 0$; the other case is similar.  Proposition \ref{app:ik}(c) gives $g_* \O_X = \O_{M_{\chi+k}}$, so the pullback map
\[ g^*\colon H^2(\O_{M_{\chi+k}}) \to H^2(\O_X) \]
is an isomorphism, so by Hodge theory the pullback map
\[ g^*\colon H^0(\Omega^2_{M_{\chi+k}}) \to H^0(\Omega^2_X) \]
is an isomorphism.  Thus if $M_{\chi+k}$ is a point then $H^0(\Omega^2_X) = 0$, and if $\dim M_{\chi+k} > 0$ then $H^0(\Omega^2_X)$ is 1-dimensional, generated by $g^* \sigma_{\chi+k}$.

It remains to show that if $\dim M_{\chi+k} > 0$ then $f^* \sigma_{\chi-k}$ is non-zero.  We have seen that $f$ is birational onto its image, so if $f^* \sigma_{\chi-k}$ vanishes then $f(X) \subset M_{\chi-k}$ is isotropic for $\sigma_{\chi-k}$, so $\dim X \le \frac12 \dim M_{\chi-k}$.  But we have seen that $\dim X = \dim M_{\chi-k} - k\chi = \dim M_{\chi+k} + k\chi$, so if $\dim M_{\chi+k} > 0$ then this is impossible.
\end{proof}

\section{Deformation condition (v)}

Like condition (iii), condition (v) will follow from an excess normal bundle calculation.  Again it is no more work, and maybe a little clearer, to prove something slightly more general:

\begin{thm}[Implies condition \ref{condv}] \label{thmv}
For each $\chi$ there is a deformation $\widetilde M_\chi$ of $M_\chi$ over $\C$, such that for each $m \ge 1$,
\begin{equation} \label{convol}
\H^*(i_{23*} E_{\chi+1} \ast i_{12*} E^{(m)}_{\chi-m}) \cong E^{(2)}_\chi[-1] \oplus E^{(2)}_\chi[m+1],
\end{equation}
where $i_{12}$ and $i_{23}$ are the closed embeddings
\[i_{12} = (1 \times i): M_{\chi-2m} \times M_\chi \to M_{\chi-2m} \times \widetilde M_\chi \]
\[i_{23} = (i \times 1): M_\chi \times M_{\chi+2} \to \widetilde M_\chi \times M_{\chi+2} \]
and $i$ is the inclusion $M_\chi \hookrightarrow \widetilde M_\chi$.
\end{thm}

\begin{proof}
We display the diagram from Proposition \ref{WZbundles}, together with the inclusion:
\[ 
\begin{tikzcd}
& W \ar[swap]{dl}{\tilde g} \ar{d}{\pi} \ar{dr}{\tilde f} & \\
X \ar[swap]{d}{e} \ar[pos=0.3]{dr}{f} & Z \ar{dl}[swap, pos=0.3]{k} \ar[pos=0.3]{dr}{\ell} & Y \ar[swap, pos=0.3]{dl}{g} \ar{d}{h} \\
M_{\chi-2m} & M_\chi \ar[hookrightarrow,swap,blue]{rd}{i} & M_{\chi+2} \\
& & {\color{blue}\widetilde M_\chi}
\end{tikzcd} \]
Because $M_{\chi+2}$ appears on the right, we see that  $\pi$ is a $\P^m$ bundle. \bigskip

As a functor, the convolution \eqref{convol} is
\[ h_*(L_Y \otimes g^* i^* i_* f_*(L_X \otimes e^*(-))). \]
In the last section we analyzed $g^* f_*$ in terms of the excess normal bundle $E$ for $W$ as the fiber product of $X$ and $Y$ over $M_\chi$.  Now we are interested in $g^* i^* i_* f_*$, so we should study the excess normal bundle $\widetilde E$ for $W$ as the fiber product of $X$ and $Y$ over $\widetilde M_\chi$.  Retracing the proof of Theorem \ref{thmiii}, we see that Theorem \ref{thmv} reduces to the assertion that
\[ \pi_*(\omega_\pi \otimes \Lambda^j \widetilde E^*[j]) = \begin{cases}
\O_Z[-m] & \text{if } j = 0 \\
\O_Z[m+1] & \text{if } j = m+1 \\
0 & \text{otherwise,}
\end{cases} \]
or equivalently (by Grothendieck--Verdier duality),
\begin{equation} \label{preclaim}
\pi_* \Lambda^j \widetilde E = \begin{cases}
\O_Z & \text{if } j = 0 \\
\O_Z[-m] & \text{if } j = m+1 \\
0 & \text{otherwise.}
\end{cases}
\end{equation}
This in turn will follow from:

\paragraph{Claim:} There is a deformation $\widetilde M_\chi$ such that $\widetilde E$ is the unique non-split extension
\begin{equation} \label{quillen_ext}
0 \to \Omega_\pi \to \widetilde E \to \O_W \to 0,
\end{equation}
sometimes called the Quillen bundle of the $\P^m$ bundle $\pi$. \bigskip

Let us argue briefly that the claim implies \eqref{preclaim}, in case it is not obvious.  To see that the extension \eqref{quillen_ext} exists and is unique, write
\[ \Ext^1(\O_W, \Omega_\pi) = H^1(\Omega_\pi); \]
then $\pi_* \Omega_\pi = \O_Z[-1]$, so $H^1(\Omega_\pi) = H^0(\O_Z) = \C$.  Next, for a point $z \in Z$ and the fiber $\pi^{-1}(z) \cong \P^m$ over it, the restriction map $H^1(\Omega_\pi) \to H^1(\Omega_{\pi^{-1}(z)})$ is an isomorphism, because the restriction map $H^0(\O_Z) \to H^0(\O_z)$ is an isomorphism.  Thus if the extension \eqref{quillen_ext} is non-split then its restriction to every fiber is non-split, so $\widetilde E|_{\pi^{-1}(z)} = \O_{\P^m}(-1)^{m+1}$, so
\[ \Lambda^j \widetilde E|_{\pi^{-1}(z)} = \O_{\P^m}(-j)^{\binom{m+1}{j}}, \]
whose cohomology vanishes for $0 < j \le m$.  On the other hand, $\Lambda^0 \widetilde E = \O_Z$ and $\Lambda^{m+1} \widetilde E = \omega_\pi$, giving \eqref{preclaim}. \bigskip

Now we prove the claim.  The excess bundle $E$ is defined by the exact sequence \eqref{excess}, and $\widetilde E$ is defined by a similar sequence, which fit together in a commutative diagram
\begin{equation} \label{comparing_excesses}
\begin{tikzcd}
0 \ar{r} & N_{W/X\times Y} \ar{r} \ar[equal]{d} & \tilde g^* f^* T_{M_\chi} \ar{r}{\rho} \ar{d} & E \ar{r} \ar{d} & 0 \\
0 \ar{r} & N_{W/X\times Y} \ar{r} & \tilde g^* f^* i^* T_{\widetilde M_\chi} \ar{r} & \widetilde E \ar{r} & 0.
\end{tikzcd}
\end{equation}
The middle vertical map $T_{M_\chi} \to i^* T_{\widetilde M_\chi}$ is injective, and its cokernel is the normal bundle of $M_\chi$ in $\widetilde M_\chi$, which is $\O_{M_\chi}$ because $\widetilde M_\chi$ is a fibration over $\C$.  The class of the extension
\[ 0 \to T_{M_\chi} \to i^* T_{\widetilde M_\chi} \to \O_{M_\chi} \to 0 \]
is the Kodaira--Spencer class of the deformation $\widetilde M$; call it $\kappa \in H^1(T_{M_\chi})$.  Applying the snake lemma to \eqref{comparing_excesses}, we find that $E$ injects into $\widetilde E$, the cokernel is $\O_W$, and in fact the class of the extension
\[ 0 \to E \to \widetilde E \to \O_W \to 0 \]
is the image of $\tilde g^* f^* \kappa$ under the map $H^1(\tilde g^* f^* T_{M_\chi}) \to H^1(E)$ induced by the map labeled $\rho$.  Thus to prove the claim we need only prove that there is a deformation $\widetilde M$, independent of $m$, such that the image of its Kodaira--Spencer class in $H^1(E)$ is non-zero. \bigskip

In Proposition \ref{excess_prop} we saw that $E = \Omega_\pi$, and that if we use the holomorphic sympletic form on $M_\chi$ to identify $T_{M_\chi}$ with $\Omega_{M_\chi}$ then the map labeled $\rho$ in \eqref{comparing_excesses} is just the restriction map on 1-forms $\tilde g^* f^* \Omega_{M_\chi} \to \Omega_\pi$.

This is surjective on $H^1$, as follows.  If we choose a point $z \in Z$, then the fiber $\pi^{-1}(z) \cong \P^m$ embeds into $M_\chi$ by Proposition \ref{WZbundles}(b).  Thus the composition
\[ H^1(\Omega_{M_\chi}) \to H^1(\tilde f^* g^* \Omega_{M_\chi}) \to H^1(\Omega_\pi) \to H^1(\Omega_{\pi^{-1}(z)}) \]
is surjective, because any K\"ahler class in $H^{1,1}(M_\chi)$ maps to a K\"ahler class in $H^{1,1}(\P^m)$, which is non-zero.  The last map $H^1(\Omega_\pi) \to H^1(\Omega_{\pi^{-1}(z)})$ is an isomorphism, as we have seen.

To conclude, choose any K\"ahler class on $M_\chi$, represent it by a hyperk\"ahler metric, let $\widetilde M_\chi \to \P^1$ be the corresponding twistor deformation, and delete the fiber over $\infty \in \P^1$ to get a deformation over $\C$.  By \cite[Prop.~25.7]{gjh}, the isomorphism $H^1(T_{M_\chi}) \cong H^1(\Omega_{M_\chi})$ induced by the holomorphic symplectic form identifies the Kodaira--Spencer class of this deformation with the K\"ahler class, up to a scalar.  This concludes the proof of the claim above, and thus of Theorem \ref{thmv}.
\end{proof}

\begin{rmk}
The scrupulous reader may worry that the twistor family $\widetilde M_\chi$ is not even K\"ahler, much less algebraic, so we must work with derived categories of coherent \emph{analytic} sheaves.  But the whole Fourier--Mukai machinery goes through undisturbed in this setting; see Ramis, Ruget, and Verdier \cite{rrv} for Grothendieck--Verdier duality, which is the hardest part.

Cautis, Kamnitzer, and Licata suggest in \cite[Rmk.~2.6]{ckl1796} that one can make do with just first-order deformation, rather than a whole deformation over $\C$, which would allow us to stay in the algebraic category.  But we could not figure out how to make this work, because the behavior of $i^* i_*$ is quite different between the inclusions $0 \hookrightarrow \C$ and $0 \hookrightarrow \Spec(\C[t]/t^2)$.
\end{rmk}


\section{Commutator condition (iv)}
\label{sectioncond4}

Condition (iv) is by far the most intricate to prove.  We recast it using Lemma \ref{weak_fours}, and assume that $\chi \ge 0$ and that $\chi-2$ lies to the right of our negative rank fix; the other cases are similar.
\begin{thm}[Implies condition \ref{condiv}]
If $\chi \geq 0$ then there is an exact triangle
\begin{equation} \label{EF-FE}
P' \to E_{\chi-1} \ast F_{\chi-1} \to F_{\chi+1} \ast E_{\chi+1},
\end{equation}
in $D(M_\chi \times M_\chi)$, where $\H^*(P') \cong \O_\Delta \otimes_\C H^*(\P^{\chi-1}).$
\end{thm}

The proof will occupy the rest of this section, and has many ingredients in common with \cite[\S6.3]{ckl1795}.

We reproduce the usual diagram
\begin{equation} \label{usu_diag_for_iv}
\begin{tikzcd}
& W \ar[swap]{ld}{\tilde g} \ar{rd}{\tilde f} \ar{d}{\pi} & \\
X\ar[swap]{d}{e} \ar{dr}[pos=0.3]{f} & Z\ar{dl} \ar{dr} & Y\ar{dl}[swap,pos=0.3]{g} \ar{d}{h} \\
M_{\chi-2} & M_{\chi} & M_{\chi+2}.
\end{tikzcd}
\end{equation}
The kernel $E_{\chi-1} \ast F_{\chi-1}$ induces the functor
\[ f_*(L_X \otimes e^* e_*(L_X^{-1} \otimes \omega_X \otimes f^*(-))), \]
so by \cite[Ex.~5.12]{huybrechts_fm} we have
\begin{equation} \label{EF}
E_{\chi-1} \ast F_{\chi-1} = (f \times f)_*[(L_X \boxtimes (L_X^{-1} \otimes \omega_X)) \otimes (e \times e)^*\O_\Delta],
\end{equation}
and similarly
\begin{equation} \label{FE}
F_{\chi+1} \ast E_{\chi+1} = (g \times g)_*[((L_Y^{-1} \otimes \omega_Y) \boxtimes L_Y) \otimes (h \times h)^*\O_\Delta].
\end{equation}
Thus it is helpful to consider the diagram
\begin{equation}
\label{wzwdiagram}
\begin{tikzcd}
& W \times_Z W \ar{dl}[swap]{\phi := \tilde g \times \tilde g} \ar{dr}{\psi := \tilde f \times \tilde f}& \\
X \times_{M_{\chi-2}} X \ar{dr}[swap]{f \times f} & & Y \times_{M_{\chi+2}} Y. \ar{dl}{g \times g}\\
& M_\chi \times M_\chi
\end{tikzcd}
\end{equation}

The outline of the argument is as follows.  The fiber product $Y \times_{M_{\chi+2}} Y$ is always irreducible of the expected dimension, and is responsible for $F_{\chi+1} \ast E_{\chi+1}$.  The fiber product $X \times_{M_{\chi-2}} X$ is responsible for $E_{\chi-1} \ast F_{\chi-1}$.  If $\chi = 0$ then it is also irreducible of the expected dimension, and is birational to $Y \times_{M_2} Y$, leading to $E_{-1} \ast F_{-1} = F_1 \ast E_1$.  If $\chi \ge 1$ then $X \times_{M_{\chi-2}} X$ has two irreducible components: the diagonal $X$ which gives rise to $P'$ in \eqref{EF-FE}, and another component $R$ which is birational to $Y \times_{M_{\chi+2}} Y$ and gives rise to $F_{\chi+1} \ast E_{\chi+1}$ in \eqref{EF-FE}.  The component $R$ always has the expected dimension.  If $\chi=1$ then the diagonal $X$ does too, which is reflected in the fact that $P' = \O_\Delta$, but if $\chi\ge 2$ then the diagonal $X$ is bigger than expected, and the wedge powers of its excess bundle lead to $P'$ being a sum of several shifted copies of $\O_\Delta$.

\subsection*{Analysis of $Y \times_{M_{\chi+2}} Y$}

\begin{lem} \label{YMY_irred}
For all $\chi \ge 0$, the fiber product $Y \times_{M_{\chi+2}} Y$ is irreducible of the expected dimension.  In particular, it is a local complete intersection, and in \eqref{FE} the derived pullback $(h \times h)^* \O_\Delta$ is just $\O_{Y \times_{M_{\chi+2}} Y}$.
\end{lem}
\begin{proof}
To lighten the notation, we set $M = M_{\chi+2}$ for the proof of this lemma, and $m = \dim M$.

The diagonal in $M \times M$ is locally cut out by a regular sequence of $m$ functions, so $Y \times_M Y$ is locally cut out of $Y \times Y$ by a sequence of $m$ functions, so the codimension of every irreducible component is at most $m$.  Moreover, because $Y \times Y$ is smooth, we see that if the codimension of every component is exactly $m$, then the sequence of $m$ functions on $Y \times Y$ is regular, so $Y \times_M Y$ is a local complete intersection, and the higher derived pullbacks $L_{>0} (h \times h)^* \O_\Delta$ vanish.

Recall the Brill--Noether stratification
\[ M = {_0 M} \supset {_1 M} \supset {_2 M} \supset \dotsb. \]
Over the locally closed set $_t M \setminus {_{t+1} M}$, we know that $Y$ is a $\P^{\chi+1+t}$ bundle.  Looking over the open set $_0 M \setminus {_1 M}$, we see that $\dim Y = m + \chi + 1$, so the expected dimension of $Y \times_M Y$ is $m + 2\chi + 2$, and the preimage of $_0 M \setminus {_1 M}$ has exactly that dimension.  We will argue that this preimage is dense in $Y \times_M Y$, because the preimage of any later $_t M \setminus {_{t+1} M}$ has smaller dimension, hence its closure cannot contribute a new irreducible component.

We have
\[ \dim {_t M} = m - t(\chi+2+t), \]
so the dimension of preimage of ${_t M} \setminus {_{t+1} M}$ in $Y \times_M Y$ is
\[ (m - t(\chi+2+t)) + 2(\chi+1+t) = (m + 2\chi + 2) - t(\chi + t). \]
When $t>0$, this is strictly less than $m + 2\chi + 2$, as desired.
\end{proof}

\begin{lem} \label{psi_is_ratl_res}
The map $\psi$ in diagram \eqref{wzwdiagram} is a rational resolution of singularities, meaning that
\begin{align*}
\psi_* \O_{W \times_Z W} &= \O_{Y \times_{M_{\chi+2}} Y} \\
\psi_* \omega_{W \times_Z W} &= \omega_{Y \times_{M_{\chi+2}} Y}.
\end{align*}
\end{lem}
\begin{proof}
First we argue that $\psi$ is a resolution of singularities.  We know that $Z$ is smooth, and $\pi\colon W \to Z$ is a $\P^1$ bundle, so $W \times_Z W$ is smooth as well. To see that $\psi$ is birational, observe that a point of $Y \times_{M_{\chi+2}} Y$ represents a sheaf $G \in M_{\chi+2}$ and a pair of 1-dimensional subspaces of $H^0(G)$:
\[ \begin{tikzcd}[row sep = tiny]
\C \ar[hookrightarrow]{dr} \\
& H^0(G), \\
\C \ar[hookrightarrow]{ur} 
\end{tikzcd} \]
whereas a point of $W \times_Z W$ represents a sheaf $G \in M_{\chi+2}$ and a pair of 1-dimensional subspaces in a 2-dimensional subspace of $H^0(G)$:
\[ \begin{tikzcd}[row sep = tiny]
\C \ar[hookrightarrow]{dr} \\
& \C^2 \ar[hookrightarrow]{r} & H^0(G). \\
\C \ar[hookrightarrow]{ur} 
\end{tikzcd} \]
The map $\psi$ forgets the 2-dimensional subspace, so it is an isomorphism away from the diagonal of $Y \times_{M_{\chi+2}} Y$: two distinct 1-dimensional subspaces are contained in a unique 2-dimensional subspace. We have seen that $Y \times_{M_{\chi+2}} Y$ is irreducible, and we find that the diagonal has codimension $\chi+1$, so in particular the complement of the diagonal is dense.

In fact we suspect that $W \times_Z W$ is the blow-up of $Y \times_{M_{\chi+2}} Y$ along the diagonal, although we do not attempt prove it. \bigskip

To prove that $\psi$ is a rational resolution, we apply \cite[Lem.~5.12]{km}.  From Lemma \ref{YMY_irred} we know that $Y \times_{M_{\chi+2}} Y$ is l.c.i.\ and hence Cohen--Macaulay, so we need to show that $R^0 \psi_* \omega_{W \times_Z W} = \omega_{Y \times_{M_{\chi+2}} Y}$.

We claim that the singular locus of $Y \times_{M_{\chi+2}} Y$ has codimension $\ge 3$, and its preimage in $W \times_Z W$ has codimension $\ge 2$.  We have just seen that $Y \times_{M_{\chi+2}} Y$ is smooth away from the diagonal, because $\psi$ is an isomorphism there.  It is also smooth away from the preimage of ${_1 M_{\chi+2}}$, that is, on the open set where $h^1(G) = 0$, because the map to $M_{\chi+2}$ is a $\P^{\chi+1} \times \P^{\chi+1}$ bundle there.  Thus the singular locus of $Y \times_{M_{\chi+2}} Y$ is contained in the preimage of $_1 M_{\chi+2}$ in the diagonal copy of $Y$, whose codimension we find to be $2\chi+3$.  Moreover we find that the preimage of $_1 M_{\chi+2}$ in the diagonal copy of $W$ has codimension $\chi+2$ in $W \times_Z W$.

Now let $V$ be the smooth locus of $Y \times_{M_{\chi+2}} Y$, let $j$ be the inclusion, let $U = \psi^{-1}(V)$, and let $i\colon U \to W \times_Z W$ be the inclusion.  Because the complements of these open sets have complements $\ge 2$, we have $R^0 i_* \omega_U = \omega_{W \times_Z W}$ and $R^0 j_* \omega_V = \omega_{Y \times_{M_{\chi+2}} Y}$ \cite[\S5]{sandor}.  Because $\psi|_U \colon U \to V$ is a proper birational morphism of smooth varieties, we have $R^0 (\psi|_U)_* \omega_U = \omega_V$ \cite[Prop.~5.77(c)]{km}.  Thus $R^0 \psi_* \omega_{W \times_Z W} = \omega_{Y \times_{M_{\chi+2}} Y}$ as desired.
\end{proof}

\begin{lem} \label{phi_vs_psi_lem}
The line bundles appearing in \eqref{EF} and \eqref{FE} satisfy
\begin{equation} \label{phi_vs_psi}
\phi^*(L_X \boxtimes (L_X^{-1} \otimes \omega_X)) = \psi^*((L_Y^{-1} \otimes \omega_Y) \boxtimes L_Y) \otimes \omega_\psi,
\end{equation}
where $\phi$ and $\psi$ are as in diagram \eqref{wzwdiagram}.
\end{lem}
\begin{proof}
First we compute $\omega_\psi$.  Because $\pi\colon W \to Z$ is smooth, we have
\[ \omega_{W \times_Z W} = \omega_\pi \boxtimes \omega_W. \]
Because $Y \times_{M_{\chi+2}} Y$ is a local complete intersection, we can speak about its canonical line bundle; because it has the expected dimension and $\omega_{M_{\chi+2}}$ is trivial, we have
\[ \omega_{Y \times_{M_{\chi+2}} Y} = \omega_Y \boxtimes \omega_Y. \]
Suppressing pullbacks from the notation for readability, we get
\[ \omega_\psi = (\omega_\pi \otimes \omega_Y^{-1}) \boxtimes (\omega_W \otimes \omega_Y^{-1}). \]
Thus the right-hand side of \eqref{phi_vs_psi} is
\[ (L_Y^{-1} \otimes \omega_\pi) \boxtimes (L_Y \otimes \omega_W \otimes \omega_Y^{-1}). \]
By construction we have $L_X \otimes L_Y = L_Z \otimes \omega_\pi$ on $W$, so this becomes
\[ (L_X \otimes L_Z^{-1}) \boxtimes (L_X^{-1} \otimes L_Z \otimes \omega_\pi \otimes \omega_W \otimes \omega_Y^{-1}). \]
Because we are on $W \times_Z W$, the $L_Z$ factor can pass through the box product and cancel to give
\begin{equation} \label{last_thing}
L_X \boxtimes (L_X^{-1} \otimes \omega_\pi \otimes \omega_W \otimes \omega_Y^{-1}).
\end{equation}
Finally, recall the exact sequence \eqref{excess}, which is
\[ 0 \to T_W \to T_X \oplus T_Y \to T_{M_\chi} \to \Omega_\pi \to 0 \]
by Proposition \ref{excess_prop}, and take determinants to get
\[ \omega_\pi \otimes \omega_X^{-1} \otimes \omega_Y^{-1} \otimes \omega_W = 0. \]
Thus \eqref{last_thing} equals the left-hand side of \eqref{phi_vs_psi}
\end{proof}

\subsection*{The case $\chi = 0$}

If $\chi = 0$ then the analogue of Lemma \ref{YMY_irred} also holds for the fiber product $X \times_{M_{-2}} X$: it is irreducible of the expected dimension, and the derived pullback $(e \times e)^* \O_\Delta$ that appears in \eqref{EF} is just $\O_{X \times_{M_{-2}} X}$.  The analogue of Lemma \ref{psi_is_ratl_res} holds as well: $\phi$ is birational, and $\phi_* \O_{W \times_Z W} = \O_{X \times_{M_{-2}} X}$.

Thus we can rewrite $E_{-1} * F_{-1}$ from \eqref{EF} as
\[ (f \times f)_* \phi_* \phi^* (L_X \boxtimes (L_X^{-1} \otimes \omega_X)) \]
using the projection formula.  By Lemma \ref{phi_vs_psi_lem} and the commutativity of diagram \eqref{wzwdiagram}, this equals
\[ (g \times g)_* \psi_* [\psi^* (L_Y \boxtimes (L_Y^{-1} \otimes \omega_Y)) \otimes \omega_\psi]. \]
By the projection formula and Lemma \ref{psi_is_ratl_res}, this equals $F_1 * E_1$ from \eqref{FE}, as desired.

\subsection*{The case $\chi = 1$}

In principle we could deal with $\chi = 1$ and $\chi \ge 2$ at the same time, but we believe that doing the simpler case $\chi = 1$ first will serve to clarify the more complicated case $\chi \ge 2$.

\begin{lem} \label{XM-1X_lci} 
If $\chi = 1$ then the fiber product $X \times_{M_{-1}} X$ has two irreducible components, both of the expected dimension: the diagonal copy of $X$, and the preimage of $_1 M_{-1}$ which we call $R$ for remainder. 

In particular, $X \times_{M_{-1}} X$ is a local complete intersection, and in \eqref{EF} the derived pullback $(e \times e)^* \O_\Delta$ is just $\O_{X \times_{M_{-1}} X}$.
\end{lem}

\begin{proof}
Because $e\colon X \to M_{-1}$ is birational, the expected dimension of the fiber product is $\dim X$.  As in the proof of Lemma \ref{YMY_irred}, every irreducible component has dimension at least that big, and the claims about l.c.i.\ and $(e \times e)^* \O_\Delta$ will follow if we show that no irreducible component is bigger.	

Over the locally closed set $_t M_{-1} \setminus {_{t+1} M_{-1}}$, we know that $X$ is a $\P^t$ bundle, so $X \times_{M_{-1}} X$ is a $\P^t \times \P^t$ bundle.  We also know that $_t M_{-1} \setminus {_{t+1} M_{-1}}$ is irreducible of dimension $\dim X - t(1+t)$, so its preimage in $X \times_{M_{-1}} X$ is irreducible of dimension $\dim X + t - t^2$.  For $t = 0$, the closure of this preimage gives one irreducible component of $X \times_{M_{-1}} X$, which is the diagonal copy of $X$.  For $t = 1$, the closure of the preimage gives a second component, which we are calling $R$.  For $t > 1$, the preimage is too small to contribute a new irreducible component.
\end{proof}

\begin{lem} \label{phi_is_ratl_res}
If $\chi = 1$ then the map $\phi$ in diagram \eqref{wzwdiagram} is birational onto $R \subset X \times_{M_{-1}} X$, and is a rational resolution of singularities:
\begin{align*}
\phi_* \O_{W \times_Z W} &= \O_R & \phi_* \omega_{W \times_Z W} &= \omega_R
\end{align*}
\end{lem}

\begin{proof}
To see that $\phi$ takes values in $R$, observe that the map $Z \to M_{-1}$ in diagram \eqref{usu_diag_for_iv} takes values in $_1 M_{-1}$.  We claim that it is birational onto $R$, and in fact is an isomorphism away from a closed set of codimension 2 in $W \times_Z W$, whose image in $R$ has codimension 3.  Observe that a point of $R$ represents a sheaf $E \in M_{-1}$ with $h^0(E) \ge 1$, hence $h^1(E) \ge 2$, together with two quotients
\[ \begin{tikzcd}[row sep = tiny]
& \C \\
H^1(E) \ar[twoheadrightarrow]{ur} \ar[twoheadrightarrow]{dr} \\
& \C,
\end{tikzcd} \]
while a point of $W \times_Z W$ represents a sheaf $E \in M_{-1}$ with quotients
\[ \begin{tikzcd}[row sep = tiny]
& & \C \\
H^1(E) \ar[twoheadrightarrow]{r} & \C^2 \ar[twoheadrightarrow]{ur} \ar[twoheadrightarrow]{dr} \\
& & \C.
\end{tikzcd} \]
The map $\phi$ forgets the $\C^2$, so it is an isomorphism away from the diagonal of $X \times_{M_{-1}} X$: perhaps it is clearest to say that two distinct 1-dimensional subspaces of $H^1(E)^*$ are contained in a unique 2-dimensional subspace.  On the other hand, away from the preimage of $_2 M_{-1}$, that is, on the open set of $R$ where $h^1(E) = 2$, there is only one possibility for the $\C^2$, so $\phi$ is an isomorphism over that open set as well.  The codimensions 2 and 3 come from analyzing the preimages of $_t M_{-1} \setminus {_{t+1} M_{-1}}$ in $Z$ and $X$ for $t \ge 2$.

Now the claim that $\phi$ is a rational resolution follows from \cite[Lem.~5.12]{km} as in the proof of Lemma \ref{psi_is_ratl_res}.  To see that $R$ is Cohen--Macaulay, apply  \cite[Thm.~21.23(b)]{eis}\footnote{In some printings of \cite{eis}, the statement of Theorem 21.23 contains a rather confusing typo: the first occurrence of $J = (0:_A I)$ should be deleted.} to $X \times_{M_{-1}} X = X \cup R$: the union is l.c.i.\ and hence Gorenstein, the first component $X$ is smooth and hence Cohen--Macaulay, so the second component $R$ is Cohen--Macaulay as well.  To see that $R^0 \phi_* \omega_{W \times_Z W} = \omega_R$, use the codimension estimates above and argue as in the proof of Lemma \ref{psi_is_ratl_res}.
\end{proof}
\smallskip

Writing $X \times_{M_{-1}} X = X \cup R$, where $X$ on the right-hand side is the diagonal, and observing that $X \cap R$ is codimension 1 in $X$, we have an exact sequence
\begin{equation} \label{mv_seq}
0 \to \O_X(-(X \cap R)) \to \O_{X \cup R} \to \O_R \to 0.
\end{equation}
If we tensor with $L_X \boxtimes (L_X^{-1} \otimes \omega_X)$ and apply $(f \times f)_*$, then the middle term becomes $E_0 * F_0$ from \eqref{EF}.  We will see that the first term becomes $\O_\Delta \in D(M_1 \times M_1)$, and the third term becomes $F_2 * E_2$ from \eqref{FE}, so the whole triangle becomes the desired \eqref{EF-FE}. \bigskip

Consider the first term of \eqref{mv_seq}.  By Proposition \ref{app:ik+1}(b), $X$ is the blow-up of $M_1$ along $_1 M_1$, and we see that $X \cap R$ is the exceptional divisor.  Because $_1 M_1$ has codimension 2 in $M_1$ and $\omega_{M_1}$ is trivial, we have $\omega_X = \O_X(X \cap R)$.  Thus when we tensor the first term of \eqref{mv_seq} with  $L_X \boxtimes (L_X^{-1} \otimes \omega_X)$, the $L_X$ and $L_X^{-1}$ factors cancel on the diagonal, and $\omega_X = \O_X(X \cap R)$ cancels with $\O_X(-(X \cap R))$, leaving $\O_X$.  Because $f\colon X \to M_1$ is a proper birational morphism of smooth varieties, or by Proposition \ref{app:ik}(c), we have $f_* \O_X = \O_{M_1}$.  Finally, we note that $f \times f$ maps the diagonal copy of $X$ onto the diagonal copy of $M_1$, so the first term of \eqref{mv_seq} becomes $\O_\Delta \in D(M_1 \times M_1)$, as desired. \bigskip

The third term of \eqref{mv_seq} is similar to what we saw with $\chi=0$: we have
\begin{align}
&(f \times f)_*[\O_R \otimes (L_X \boxtimes (L_X^{-1} \otimes \omega_X))] \notag \\
&= (f \times f)_* \phi_* \phi^* (L_X \boxtimes (L_X^{-1} \otimes \omega_X)) \label{proj_fmla_argument} \\
&= (g \times g)_* \psi_*[\psi^* (L_Y \boxtimes (L_Y^{-1} \otimes \omega_Y)) \otimes \omega_\psi], \notag
\end{align}
where in the second line we have used Lemma \ref{phi_is_ratl_res} and the projection formula, and in the third we have used Lemma \ref{phi_vs_psi_lem} and the commutativity of diagram \eqref{FE}.  Now by the projection formula and Lemma \ref{psi_is_ratl_res}, this equals $F_2 * E_2$ from \eqref{FE}, as desired.

\subsection*{The case $\chi \ge 2$}

Now the map $e\colon X \to M_{\chi-2}$ is birational onto its image $_1 M_{\chi-2}$, which has codimension $\chi-1$ in $M_{\chi-2}$, so the expected dimension of the fiber product $X \times_{M_{\chi-2}} X$ is $\dim X - \chi + 1$.  Emulating the proof of Lemma \ref{XM-1X_lci} we find that $X \times_{M_{\chi-2}} X$ has two irreducible components: the diagonal $X$, which is bigger than expected, and the preimage of $_2 M_{\chi-2}$, which has the expected dimension and which we again call $R$ for remainder.

Because $X \times_{M_{\chi-2}} X$ has two components of different dimensions, we do not know how to deal with the term $(e \times e)^* \O_\Delta$ in \eqref{EF} directly.  We work around this by blowing up one copy of $X$.  We augment diagram \eqref{usu_diag_for_iv} as follows:
\begin{equation}
\label{extendeddiagram}
\begin{tikzcd}[column sep = 1.6 em]
& & & B \ar[swap]{dl}{\beta} \ar{drr}{b} \ar{d}{\rho} & & & W \ar{dl} \ar{dr} \ar{d} \\
& & A\ar{dll} \ar{drr}[swap]{a} & C \ar{dlll} \ar{drrr}[pos=0.3]{c} & & X \ar{dl}[swap]{e} \ar{dr}[pos=0.6]{f} & Z \ar{dll} \ar{drr} & Y \ar{dl} \ar{dr}\\
M_{-\chi} & & & & M_{\chi-2} & & M_\chi & & M_{\chi+2}.
\end{tikzcd}
\end{equation}
Here $A$ and $C$ are the usual correspondences and $B$ is the fiber product $A \times_{M_{\chi-2}} X$.  The map $\rho$ is a $\P^{\chi-1}$ bundle by Proposition \ref{WZbundles}, and is in some sense a resolution of $f$ which is only generically a $\P^{\chi-1}$ bundle.  More precisely, we know that $C$ is the blow-up of $M_\chi$ along $_1 M_\chi$ by Proposition \ref{app:ik+1}(b), and we claim that $\rho$ is the proper transform of $f$.  To see this, note that $B$ is not quite the fiber product $C \times_{M_\chi} X$, but it is the irreducible component of that fiber product containing the preimage of $M_\chi \setminus {_1 M_\chi}$.  In particular, $B$ is the blow-up of $X$ along $f^{-1}({_1 M_\chi})$.

We augment diagram \eqref{wzwdiagram} as follows:
\begin{equation}
\label{bigdiagram}
\begin{tikzcd}[column sep=0.5]
\color{blue} B \cup R' \ar[blue]{d} \ar[hookrightarrow,blue]{r} \ar[blue]{drr} & B \times X \ar{d}{\beta \times 1} \ar[sloped]{drrr}{b \times 1} & & & W \times_Z W \ar{dll} \ar{dr} &\\
\color{blue} B \ar[hookrightarrow,blue]{r} \ar[blue]{d} & A \times X \ar{d}{a \times e} & \color{blue} X \cup R \ar[hookrightarrow,blue]{rr} \ar[blue]{dll} \ar{drr} & & X \times X \ar[sloped,swap]{dlll}{e \times e} \ar{d}{f \times f} & Y \times_{M_{\chi+2}} Y \ar{dl} \\
\color{blue} \Delta \ar[hookrightarrow,blue]{r} & M_{\chi-2} \times M_{\chi-2} & & & M_\chi \times M_\chi
\end{tikzcd}
\end{equation}
The terms highlighted in blue are the diagonal in $M_{\chi-2} \times M_{\chi-2}$ and its preimages in $X \times X$, $A \times X$, and $B \times X$.  In the top left we are asserting:

\begin{lem} \label{BMX_lci} 
The fiber product $B \times_{M_{\chi-2}} X$ has two irreducible components: the graph of $b$ inside $B \times X$, which is isomorphic to $B$ and is the preimage of the diagonal $X \subset X \times X$; and another component $R'$ which is the preimage of $R \subset X \times X$.  As the preimage of $\Delta \subset M_{\chi-2} \times M_{\chi-2}$, both have greater than expected dimension, but as the preimage of $B \subset A \times X$, both have the expected dimension.

In particular, $B \times_{M_{\chi-2}} X$ is a local complete intersection, and the derived pullback $(\beta \times 1)^* \O_B$ is just $\O_{B \times_{M_{\chi-2}} X}$.  More generally, any vector bundle supported on $B \subset A \times X$ has no higher derived pullbacks via $\beta \times 1$.
\end{lem}
\begin{proof}
Because $B$ is smooth and $\beta$ is birational, every irreducible component of $(\beta \times 1)^{-1}(B)$ has dimension greater than or equal to $\dim B$.  As in the proof of Lemma \ref{XM-1X_lci}, for each stratum $_t M_{\chi-2} \setminus {_{t+1} M_{\chi-2}}$ in the diagonal of $M_{\chi-2} \times M_{\chi-2}$, we can compute the dimension of its preimage in $B \times X$.  We find that $t=1$ and $t=2$ give components of the same dimension as $B$, while for $t \ge 3$ the preimages are too small to contribute a new component.  Then the remaining claims follow as before.
\end{proof}

Because $b\colon B \to X$ is a proper birational morphism of smooth varieties, we have $b_* \O_B = \O_X$.  Thus in computing \eqref{EF} we have
\begin{align*}
(e \times e)^* \O_\Delta
&= (b \times 1)_* (b \times 1)^* (e \times e)^* \O_\Delta \\
&= (b \times 1)_* (\beta \times 1)^* (a \times e)^* \O_\Delta.
\end{align*}
In \S\ref{sectioncond3} we saw that
\[ \H^q((a \times e)^* \O_\Delta) = \Lambda^{-q} T_\rho, \]
supported on $B \subset A \times X$.  Thus
\[ \H^q((\beta \times 1)^* (a \times e)^* \O_\Delta) = (\beta \times 1)^* \Lambda^{-q} T_\rho, \]
supported on $B \cup R' \subset B \times X$, using the last part of Lemma \ref{XM-1X_lci}.

Consider the exact sequence
\[ 0 \to \O_B(-B \cap R') \to \O_{B \cup R'} \to \O_{R'} \to 0. \]
Suppose we tensor with $(\beta \times 1)^* \Lambda^{-q} T_\rho$, so the middle term becomes $\H^q((\beta \times 1)^* (a \times e)^* \O_\Delta)$, and we apply the functor that turns $(\beta \times 1)^* (a \times e)^* \O_\Delta$ into $E_{\chi-1} * F_{\chi-1}$: that is, apply $(b \times 1)_*$, then tensor with $L_X \boxtimes (L_X^{-1} \otimes \omega_X)$, then apply $(f \times f)_*$.  In Lemmas \ref{first_term_lemma} and \ref{third_term_lemma} below, we will see that the first term becomes $\O_\Delta[-\chi+1-q] \in D(M_\chi \times M_\chi)$, and the third becomes $F_{\chi+1} * E_{\chi+1}$ if $q=0$, and zero if $q < 0$.  We would like to say that this gives the desired exact triangle \eqref{EF-FE}, but this is not quite right.

Instead, we should consider the map from $(\beta \times 1)^* (a \times e)^* \O_\Delta$ to its zeroth cohomology sheaf $\O_{B \cup R'}$ and then onto $\O_{R'}$, and let $K$ be shifted cone defined by the exact triangle
\begin{equation} \label{clever_triangle}
K \to (\beta \times 1)^* (a \times e)^* \O_\Delta \to \O_{R'}
\end{equation}
in $D(B \times X)$.  If we apply $(b \times 1)_*$, tensor with $L_X \boxtimes (L_X^{-1} \otimes \omega_X)$, and apply $(f \times f)_*$, then the middle term becomes $E_{\chi-1} * F_{\chi-1}$ from \eqref{EF}.  The third term becomes $F_{\chi+1} * E_{\chi+1}$: using Lemma \ref{third_term_lemma}(a) below, we write
\begin{align*}
&(f \times f)_*[(b \times 1)_* \O_{R'} \otimes (L_X \boxtimes (L_X^{-1} \otimes \omega_X))] \\
&=(f \times f)_*[\phi_* \O_{W \times_Z W} \otimes (L_X \boxtimes (L_X^{-1} \otimes \omega_X))],
\end{align*}
then we proceed as in \eqref{proj_fmla_argument} to get $F_{\chi+1} * E_{\chi+1}$.  For the first term, we have $\H^0(K) = \O_B(-B \cap R')$, which becomes $\O_\Delta[-\chi+1]$ by Lemma \ref{first_term_lemma} below.  For $q < 0$, we have $\H^q(K) = (\beta \times 1)^* \Lambda^{-q} T_\rho$, which fits into an exact sequence
\[ 0 \to \O_B(-B \cap R') \otimes \Lambda^{-q} T_\rho \to \H^q(K) \to \O_{R'} \otimes (\beta \times 1)^* \Lambda^{-q} T_\rho \to 0. \]
By Lemma \ref{first_term_lemma} below, the first term becomes $\O_\Delta[-\chi+1-q]$, and the third term vanishes.  Thus a Grothendieck spectral sequence like the one used in \S\ref{sectioncond3} shows that $K$ becomes an object $P'$ with $\H^*(P') = \O_\Delta \otimes \H^*(\P^{\chi-1})$, and the triangle \eqref{clever_triangle} becomes \eqref{EF-FE}, as desired.

It remains to prove the two promised lemmas.
\pagebreak 

\begin{lem} \label{first_term_lemma}
\begin{multline*}
(f \times f)_* [(b \times 1)_*(\O_B(-B \cap R') \otimes \Lambda^{-q} T_\rho) \otimes (L_X \boxtimes (L_X^{-1} \otimes \omega_X))] \\
= \O_\Delta[-\chi+1-q] \in D(M_\chi \times M_\chi).
\end{multline*}
\end{lem}
\begin{proof}
Because $B \subset B \times X$ lies over the diagonal of $X \times X$, the factors $L_X$ and $L_X^{-1}$ again cancel, and we can rewrite the whole expression as
\begin{equation} \label{thing_to_simplify}
\Delta_* c_* \rho_*(\O_B(-B \cap R') \otimes \Lambda^{-q} T_\rho \otimes b^* \omega_X)
\end{equation}

We claim that $b^* \omega_X \otimes O_B(-B \cap R') = \omega_\rho$.  Looking at diagram \eqref{extendeddiagram}, we have seen that $C$ is the blow-up of $M_\chi$ along $_1 M_\chi$, which has codimension $\chi+1$, so if we let $E := e^{-1}({_1 M_\chi})$ denote the exceptional divisor, then $\omega_C = \O_C(\chi \cdot E)$.  We have also seen that $B$ is the blow-up of $X$ along $f^{-1}({_1 M_\chi})$, which has codimension $\chi$ in $X$, so $\omega_B = b^* \omega_X \otimes \rho^* \O_C((\chi-1) E)$.  Thus $\omega_\rho = b^* \omega_X \otimes \rho^* \O_C(-E)$.  We also have $\rho^{-1}(E) = B \cap R'$, so the claim is proved.

Now \eqref{thing_to_simplify} becomes
\[ \Delta_* c_* \rho_* \Lambda^{\chi-1+q} \Omega_\rho. \]
Pushing down to $C$ we get $\O_C[-\chi+1-q]$.  Because $c$ is a proper birational morphism of smooth varieties, or by Proposition \ref{app:ik}(c), we have $c_* \O_C = \O_{M_\chi}$, and the lemma follows.
\end{proof}

\begin{lem} \label{third_term_lemma} \ 
\begin{enumerate}
\renewcommand \labelenumi {(\alph{enumi})}
\item $(b \times 1)_* \O_{R'} = \phi_* \O_{W \times_Z W}$.
\item $(b \times 1)_*(\O_{R'} \otimes (\beta \times 1)^* \Lambda^j T_\rho) = 0$ for $j > 0$.
\end{enumerate}
\end{lem}

\begin{proof}
We reluctantly consider a diagram
\[ \begin{tikzcd}[column sep={3em,between origins}]
& & \Upsilon' \ar[swap]{ld}{\upsilon'} \ar{rrrd}{u'} \\
& \Upsilon \ar[swap]{ld}{\upsilon} \ar{rrrd}{u} & & & & \Upsilon'' \ar[swap]{ld}{\upsilon''} \\
\color{blue} B \cup R' \ar[blue]{rrrd}{b \times 1} & & & & W \times_Z W \ar[swap]{ld}{\phi} \\
& & & \color{blue} X \cup R
\end{tikzcd} \]
which will fit on the top of diagram \eqref{bigdiagram}. 

A point of $\Upsilon$ represents a sheaf $E \in M_{\chi-2}$ with quotients
\[ \begin{tikzcd}[row sep = tiny]
& & \C \\
H^1(E) \ar[twoheadrightarrow]{r} & \C^2 \ar[twoheadrightarrow]{ur} \ar[twoheadrightarrow]{dr} \\
& & \C,
\end{tikzcd} \]
and a subspace $\C^{\chi-1} \subset H^0(E)$.  A point of $\Upsilon'$ represents a sheaf $E$ with quotients as above and a flag of subspaces $\C^{\chi-1} \subset \C^\chi \subset H^0(E)$.  A point of $\Upsilon''$ represents a sheaf $E$ with quotients as above and a subspace $\C^\chi \subset H^0(E)$.  All three spaces are smooth, being bundles over fiber products that are smooth by Proposition \ref{WZbundles}(a).

The map $\upsilon\colon \Upsilon \to R'$ that forgets the $\C^2$ is a rational resolution of singularities: just as in the proof of Lemma \ref{phi_is_ratl_res}, we find that $R'$ is Cohen--Macaulay because $B \cup R'$ is l.c.i.\ and $B$ is smooth, and we find that $\upsilon$ is an isomorphism away from a set of codimension 3 in $\Upsilon$ and codimension 4 in $R'$, so $\upsilon_* \omega_\Upsilon = \omega_{R'}$.

The map $\upsilon'\colon \Upsilon' \to \Upsilon$ that forgets the $\C^\chi$ is a proper birational morphism of smooth varieties, so $\upsilon'_* \O_{\Upsilon'} = \O_\Upsilon$.  Similarly, $\upsilon''_* \O_{\Upsilon''} = \O_{W \times_Z W}$.

The map $u'\colon \Upsilon' \to \Upsilon''$ that forgets the $\C^{\chi-1}$ is a $\P^{\chi-1}$ bundle, so $u'_* \O_{\Upsilon'} = \O_{\Upsilon''}$.  Thus part (a) of the lemma is proved by pushing down $\O_{\Upsilon'}$ in two different ways.  Part (b) will follow once we show that $u'_* \Lambda^j T_\rho = 0$ for $j > 0$.\bigskip

Let $y \in \Upsilon''$ be a point representing a sheaf $E \in M_{\chi-2}$ with quotients as above and subspace $\C^\chi \subset H^0(E)$.  The fiber $u'^{-1}(y)$ is a $\Gr(\chi-1,\chi)$, parametrizing subspaces $W \subset \C^\chi \subset H^0(E)$.  To understand how this fiber maps to $B$ and how it relates to the fiber of $\rho\colon B \to C$, recall that the first 1-dimensional quotient of $E$ gives an extension
\[ 0 \to \O_S \to F \to E \to 0 \]
with $F \in M_\chi$, and 1-dimensional subspace $\ell \subset H^0(F)$.  The flag $W \subset \C^\chi \subset H^0(E)$ gives a flag $\ell \subset W' \subset \C^{\chi+1} \subset H^0(F)$, where $W' \cong \C^\chi$.  The map $\rho\colon B \to C$ forgets the 1-dimensional subspace $\ell$ but remembers $W'$, so the tangent space to the fiber of $\rho$ is $\Hom(\ell, W'/\ell) = \Hom(\ell, W)$.  But back in the fiber of $\upsilon'$, the space $W$ is varying while $\ell$ is fixed, so the pullback of $T_\rho$ is just the tautological sub-bundle of the Grassmannian $\Gr(\chi-1,\chi)$.  The exterior powers of this have no cohomology by the Borel--Weil--Bott theorem, so $u'_* \Lambda^j T_\rho = 0$ for $j > 0$, as desired.
\end{proof}

\appendix

\section{On Grassmannians of coherent sheaves} \label{appendix}

Let $X$ be an integral, 
Cohen--Macaulay, 
quasi-projective variety 
over $\C$. 
Let $F$ be a coherent sheaf of rank $r > 0$, possibly twisted by a Brauer class $\alpha \in \Br(X)$.  For foundations of twisted sheaves we refer to C\u{a}ldar\u{a}ru \cite{andrei}.

Suppose that $F$ has projective dimension at most 1, and fix a two-step resolution by ($\alpha$-twisted) vector bundles
\begin{equation} \label{two_step_res}
0 \to V_0 \xrightarrow\rho V_1 \to F \to 0,
\end{equation}
with $\rank(V_0) = s$ and thus $\rank(V_1) = r+s$.

Consider the stratification
\[ X = X_0 \supset X_1 \supset X_2 \supset \dotsb, \]
where
\[ X_i = \{ x \in X : \dim(F|_x) \ge r+i \} \]
with a scheme structure coming from the $(r+i-1)^\th$ Fitting ideal of $F$.  From \cite[Ex.~10.9]{eis}, we see that if $X_i$ is non-empty then its codimension is at most $i(r+i)$.  If the codimension of $X_i$ is exactly this ``expected codimension,'' then $F$ satisfies the hypotheses of every result below.

For a fixed integer $k$ with $0 < k \le r$, consider the Grassmannian of rank $k$ quotients
\[ \Gr(F,k) \xrightarrow{\ \pi\ } X, \]
characterized by the universal property that a map $T \to \Gr(F,k)$ is the same as a map $f\colon T \to X$ and a surjection $f^* F \twoheadrightarrow E$ onto an $f^* \alpha$-twisted vector bundle $E$ of rank $k$.  This implies the fibers of $\pi$ are Grassmannians of the fibers of the sheaf $F$: that is, for any point $x \in X$ we have $\pi^{-1}(x) = \Gr(F|_x, k)$.  For foundations of Grassmannians of (quasi-)coherent sheaves we refer to \cite[\S9.7]{egas}.\footnote{Note that this is the 1971 edition of EGA I.  The more widely available 1960 edition does not cover Grassmannians of quasi-coherent sheaves, only projectivizations.}

The case $k=1$ is the most familiar:
\[ \Gr(F,1) = \P F := \relProj(\Sym^*\!F). \]
This is discussed, for example, in \cite[Thm.~III-44]{eh} or \cite[\href{https://stacks.math.columbia.edu/tag/01O4}{Tag 01O4}]{stacks-project}.
\newpage 

\begin{prop} \label{app:ik}
Suppose that $\codim(X_i) \ge ik$ for all $i > 0$.  Then:
\begin{enumerate}
\renewcommand \labelenumi {(\alph{enumi})}
\item $\pi$ is a local complete intersection morphism.
\item The relative canonical bundle $\omega_\pi$ is given by
\[ \pi^*(\det F)^k \otimes (\det Q)^{-r}, \]
where $Q$ is the tautological quotient bundle on $\Gr(F,k)$.
\item $\pi_* \O_{\Gr(F,k)} = \O_X$.
\end{enumerate}
\end{prop}
\begin{proof}
In part (b), note that $Q$ is a $\pi^* \alpha$-twisted vector bundle, so $\det Q$ is twisted by $\pi^* \alpha^k$, and $\det F$ is twisted by $\alpha^r$; thus our formula for $\omega_\pi$ gives a naturally untwisted line bundle. \bigskip

(a) Consider
\[ \Gr(V_1,k) \xrightarrow{\ \varpi\ } X, \]
and let $Q'$ be the ($\varpi^* \alpha$-twisted) tautological quotient bundle on $\Gr(V_1,k)$.  Then $\Gr(F,k)$ is a closed subscheme of $\Gr(V_1,k)$ by \cite[Prop~9.7.9]{egas}, and the proof of [ibid., Lem.~9.7.9.1] can be reinterpreted as saying that it is cut out by the section of $\sheafHom(\varpi^* V_0, Q')$ given by the composition
\[ \varpi^* V_0 \xrightarrow{\varpi^* \rho} \varpi^* V_1 \twoheadrightarrow Q'. \]
The idea is that $\varpi^* V_1 \twoheadrightarrow Q'$ factors through $\varpi^* F$ if and only if the composite $\varpi^* V_0 \to Q'$ is zero.

We claim that this section of the vector bundle $\sheafHom(\varpi^* V_0, Q')$ is regular.  Because $X$ is Cohen--Macaulay and $\varpi$ is smooth, $\Gr(V_1,k)$ is Cohen--Macaulay, so it is enough to show that the codimension of $\Gr(F,k)$ equals the rank of $\sheafHom(\varpi^* V_0, Q')$, that is, $ks$.  Clearly every component has codimension at most $ks$.  We have $\dim \Gr(V_1,k) = \dim X + k(r+s-k)$, and the preimage of $X \setminus X_1$ in $\Gr(F,k)$ has dimension $\dim X + k(r-k)$, so it has codimension $ks$ in $\Gr(V_1,k)$.  For $i \ge 1$, the preimage of $X_i \setminus X_{i+1}$ in $\Gr(F,k)$ has dimension
\begin{align*}
&\dim(X_i \setminus X_{i+1}) + k(r+i-k) \\
&\le \dim X - ik + k(r+i-k) \\
&= \dim X + k(r-k),
\end{align*}
so either it contributes a new irreducible component of codimension $ks$, or it is too small to contribute a new irreducible component. \bigskip

(b) We use the adjunction formula
\[ \omega_\pi = \omega_\varpi|_{\Gr(F,k)} \otimes \det N, \]
where $N$ is the normal bundle of $\Gr(F,k)$ in $\Gr(V_1,k)$.  The relative canonical bundle of $\varpi$ is well-known to be
\[ \omega_\varpi = \varpi^*(\det V_1)^k \otimes (\det Q')^{-r-s}. \]
The restriction of $Q'$ to $\Gr(F,k)$ is $Q$, so
\[ \omega_\varpi|_{\Gr(F,k)} = \pi^*(\det V_1)^k \otimes (\det Q)^{-r-s}. \]
Because $\Gr(F,k)$ is cut out by a regular section, we have
\[ N = (\varpi^* V_0^\vee \otimes Q')|_{\Gr(F,k)} = \pi^* V_0^\vee \otimes Q', \]
so
\[ \det N = \pi^*(\det V_0)^{-k} \otimes (\det Q)^s. \]
Now the proposition follows from $\det F = \det V_1 \otimes (\det V_0)^{-1}$. \bigskip

(c) We have seen that $\Gr(F,k)$ is cut out of $\Gr(V_1,k)$ by a regular section of $\varpi^* V_0^\vee \otimes Q'$, so we get a Koszul resolution
\[ 0 \to \Lambda^{ks}(\varpi^* V_0 \otimes Q'^\vee) \to \dotsb \to \O_{\Gr(V_1,k)} \to \O_{\Gr(F,k)} \to 0. \]
We know that $\varpi_* \O_{\Gr(V_1,k)} = \O_X$.  We claim that $\varpi_* \Lambda^l (\varpi^* V_0 \otimes Q'^\vee) = 0$ for $1 \le l \le ks$.

By \cite[Lem.~0.5]{kapranov_gr}, we have
\[ \Lambda^l(\varpi^* V_0 \otimes Q'^\vee) = \bigoplus_{|\alpha|=l} \varpi^* \Sigma^\alpha V_0 \otimes \Sigma^{\alpha^*} Q'^\vee, \]
where the direct sum is over Young diagrams $\alpha$ with $l$ boxes, $\alpha^*$ denotes the transposed Young diagram, and $\Sigma^\alpha$ and $\Sigma^{\alpha^*}$ are Schur functors.  If $\Sigma^\alpha V_0$ is not zero then $\alpha$ can have at most $s$ rows, so $\alpha^*$ can have at most $s$ columns.  We have assumed that $k \le r$, so $s \le r+s-k$, and it follows from Kapranov's Borel--Weil--Bott calculation \cite[Prop.~2.2(b)]{kapranov_gr} that if $|\alpha| \ge 1$ and $\alpha^*$ has at most $r+s-k$ columns then $\varpi_* \Sigma^{\alpha^*} Q'^\vee = 0$.  By the projection formula, this implies the claim.
\end{proof}
\pagebreak 

\begin{prop} \label{app:ik+1}
Suppose that $\codim(X_i) \ge ik+1$ for all $i > 0$.  Then:
\begin{enumerate}
\renewcommand \labelenumi {(\alph{enumi})}
\item $\Gr(F,r)$ is irreducible.
\item If $k=r$, then $\Gr(F,r)$ is the blow-up of $X$ along $X_1$.
\end{enumerate}
\end{prop}
\begin{proof}
(a) Retracing the proof of Proposition \ref{app:ik}(a), we see that for $i > 0$, the preimage of $X_i \setminus X_{i+1}$ in $\Gr(F,k)$ is too small to contribute a new irreducible component. \bigskip

(b) Let $I_{X_1}$ denote the ideal sheaf of $X_1$.  Because $X_1$ has the expected codimension, the Eagon--Northcott complex of the map $\rho^\top\colon V_1^\vee \to V_0^\vee$ is exact:
\[ 
\dotsb \to \Sym^2 V_0 \otimes \Lambda^{s+2} V_1^\vee
\to V_0 \otimes \Lambda^{s+1} V_1^\vee
\to \Lambda^s V_1^\vee
\to \Lambda^s V_0^\vee \otimes I_{X_1} \to 0.
\]
Tensoring with $\det V_1 = \Lambda^{r+s} V_1$, we get an exact sequence
\[ 
\dotsb \to \Sym^2 V_0 \otimes \Lambda^{r-2} V_1
\to V_0 \otimes \Lambda^{r-1} V_1
\to \Lambda^r V_1
\to \det F \otimes I_{X_1} \to 0,
\]
where the boundary maps are very easy to describe: a basic element
\[ (b_1 \dotsb b_d) \otimes (a_{d+1} \wedge \dotsb \wedge a_r) \in \Sym^d V_0 \otimes \Lambda^{r-d} V_1 \]
goes to
\[ \sum_{i=1}^d (b_1 \dotsb \widehat{b_i} \dotsb b_d) \otimes (\rho(b_i) \wedge a_{d+1} \wedge \dotsb \wedge a_r) \in \Sym^{d-1} V_0 \otimes \Lambda^{r-d+1} V_1. \]
This gives an embedding of $\P(I_{X_1} \otimes \det F) = \P(I_{X_1})$ into $\P \Lambda^r V_1$.

We claim that this $\P(I_{X_1})$ contains $\Gr(F,r) \subset \Gr(V_1,r)$ in its Pl\"ucker embedding in $\P \Lambda^r V_1$.  We have seen that $\Gr(F,r)$ is cut out of $\Gr(V_1,r)$ by a section of $\varpi^* V_0^\vee \otimes Q'$.  Similarly, $\P(I_{X_1})$ is cut out of $\P \Lambda^r V_1$ by a section of $\sheafHom(V_0 \otimes \Lambda^{r-1} V_1, L)$, where $L$ is the universal quotient line bundle on $\P \Lambda^r V_1$.  We have $\det Q' = L|_{\Gr(V_1,r)}$.  If the first section vanishes, that is, if $V_1 \to Q'$ annihilates the image of $V_0 \to V_1$, then also $\Lambda^r V_1 \to \Lambda^r Q'$ annihilates the image of $V_0 \otimes \Lambda^{r-1} V_1 \to \Lambda^r V_1$, so the second section vanishes as well.

Now we know that
\[ \Bl_{X_1}(X) = \relProj(\O_X \oplus I_{X_1} \oplus I_{X_1}^2 \oplus \dotsb) \]
is the irreducible component of 
\[ \P(I_{X_1}) = \relProj(\O_X \oplus I_{X_1} \oplus \Sym^2 I_{X_1} \oplus \dotsb) \]
that contains the preimage of $X \setminus X_1$.  But $\Gr(F,r)$ contains this preimage, and we have seen that it is irreducible, so it coincides with $\Bl_{X_1}(X)$.
\end{proof}

\begin{prop} \label{app:ik+2}
Suppose that $X$ and $\Gr(F,k)$ are smooth.
\begin{enumerate}
\renewcommand \labelenumi {(\alph{enumi})}
\item Suppose that $k < r$.  If $\codim(X_i) \ge ik+2$ for all $i > 0$, then $\Pic(\Gr(F,k)) = \Pic(X) \oplus \Z$.
\item Suppose that $k = r$.  If $\codim(X_1) = r+1$ as expected, $\codim(X_i) \ge ik+2$ for all $i > 1$, and $X_1\setminus X_2$ is smooth, then the same conclusion holds.
\end{enumerate}
\end{prop}
\begin{proof}
(a) For $i \ge 1$, the stratum $X_i$ has codimension at least $ik+2$ in $X$ by hypothesis, and calculating as in the proof of Proposition \ref{app:ik}(a) we find that $\pi^{-1}(X_i \setminus X_{i+1})$ has codimension at least 2 in $\Gr(F,k)$.  Thus the restrictions
\[ \Pic(X) \to \Pic(X \setminus X_1) \qquad \text{and} \qquad \Pic(\Gr(F,k)) \to \Pic(\pi^{-1}(X \setminus X_1)) \]
are isomorphisms.  Now $\pi^{-1}(X \setminus X_1) \to X \setminus X_1$ is a Grassmannian bundle over a smooth base, so we have an exact sequence
\[ 0 \to \Pic(X \setminus X_1) \xrightarrow{\pi^*} \Pic(\pi^{-1}(X \setminus X_1)) \to \Pic(\Gr(r,k)) \to \Br(X \setminus X_1), \]
where the second map restricts a line bundle on $X \setminus X_1$ to a fiber, and the third map sends the generator of $\Pic(\Gr(r,k)) = \Z$ to the Brauer class $\alpha|_{X \setminus X_1}$.  The kernel of the latter map is a finite-index subgroup of $\Z$, that is, another copy of $\Z$, and so we have an exact sequence
\[ 0 \to \Pic(X \setminus X_1) \xrightarrow{\pi^*} \Pic(\pi^{-1}(X \setminus X_1)) \to \Z \to 0, \]
which necessarily splits. \bigskip

(b) Now if $i \ge 2$ then $X_i \subset X$ and $\pi^{-1}(X_i) \subset \Gr(F,k)$ both have codimension $\ge 2$, so the restrictions
\[ \Pic(X) \to \Pic(X \setminus X_2) \qquad \text{and} \qquad \Pic(\Gr(F,r)) \to \Pic(\pi^{-1}(X \setminus X_2)) \]
are isomorphisms.  By Proposition \ref{app:ik+1}(b), $\pi^{-1}(X \setminus X_2)$ is the blow-up of $X \setminus X_2$ along the smooth center $X_1 \setminus X_2$, so again
\[ \Pic(\pi^{-1}(X \setminus X_2)) = \Pic(X \setminus X_2) \oplus \Z. \qedhere \]
\end{proof}

\let\bloot\section 
\def\section*#1{\bloot*{#1} \addcontentsline{toc}{section}{#1}} 
\bibliographystyle{plain}
\bibliography{main}

\newcommand \httpurl [1] {\href{https://#1}{\nolinkurl{#1}}}
\begin{thebibliography}{10}

\bibitem{adm}
N.~Addington, W.~Donovan, and C.~Meachan.
\newblock Moduli spaces of torsion sheaves on {K}3 surfaces and derived
  equivalences.
\newblock {\em J. Lond. Math. Soc. (2)}, 93(3):846--865, 2016.
\newblock Also \href{https://arxiv.org/abs/1507.02597}{arXiv:1507.02597}.

\bibitem{baranovsky}
V.~Baranovsky.
\newblock Moduli of sheaves on surfaces and action of the oscillator algebra.
\newblock {\em J. Differential Geom.}, 55(2):193--227, 2000.
\newblock Also \href{https://arxiv.org/abs/math/9811092}{math/9811092}.

\bibitem{beauville_counting}
A.~Beauville.
\newblock Counting rational curves on {K3} surfaces.
\newblock {\em Duke Math. J.}, 97(1):99--108, 1999.
\newblock Also \href{https://arxiv.org/abs/alg-geom/9701019}{alg-geom/9701019}.

\bibitem{ckl1795}
S.~Cautis, J.~Kamnitzer, and A.~Licata.
\newblock Categorical geometric skew {H}owe duality.
\newblock {\em Invent. Math.}, 180(1):111--159, 2010.
\newblock Also \href{https://arxiv.org/abs/0902.1795}{arXiv:0902.1795}.

\bibitem{ckl1796}
S.~Cautis, J.~Kamnitzer, and A.~Licata.
\newblock Coherent sheaves and categorical {$\mathfrak{sl}_2$} actions.
\newblock {\em Duke Math. J.}, 154(1):135--179, 2010.
\newblock Also \href{https://arxiv.org/abs/0902.1796}{arXiv:0902.1796}.

\bibitem{ckl1797}
S.~Cautis, J.~Kamnitzer, and A.~Licata.
\newblock Derived equivalences for cotangent bundles of {G}rassmannians via
  categorical {$\mathfrak{sl}_2$} actions.
\newblock {\em J. reine angew. Math.}, 675:53--99, 2013.
\newblock Also \href{https://arxiv.org/abs/0902.1797}{arXiv:0902.1797}.

\bibitem{ckl_quiver}
S.~Cautis, J.~Kamnitzer, and A.~Licata.
\newblock \vphantom{Z}
  {C}oherent sheaves on quiver varieties and categorification.
\newblock {\em Math. Ann.}, 357(3):805--854, 2013.
\newblock Also \href{https://arxiv.org/abs/1104.0352}{arXiv:1104.0352}.

\bibitem{andrei}
A.~C\u{a}ldar\u{a}ru.
\newblock {\em Derived categories of twisted sheaves on {C}alabi--{Y}au
  manifolds}.
\newblock PhD thesis, Cornell, 2000.
\newblock Available at
  \httpurl{people.math.wisc.edu/~andreic/publications/ThesisSingleSpaced.pdf}.

\bibitem{eis}
D.~Eisenbud.
\newblock {\em Commutative algebra}, volume 150 of {\em Graduate Texts in
  Mathematics}.
\newblock Springer-Verlag, New York, 1995.
\newblock With a view toward algebraic geometry.

\bibitem{eh}
D.~Eisenbud and J.~Harris.
\newblock {\em The geometry of schemes}, volume 197 of {\em Graduate Texts in
  Mathematics}.
\newblock Springer-Verlag, New York, 2000.

\bibitem{fulton}
W.~Fulton.
\newblock {\em Intersection theory}, volume~2 of {\em Ergebnisse der Mathematik
  und ihrer Grenzgebiete. 3. Folge.}
\newblock Springer-Verlag, Berlin, second edition, 1998.

\bibitem{grojnowski}
I.~Grojnowski.
\newblock Instantons and affine algebras. {I}. {T}he {H}ilbert scheme and
  vertex operators.
\newblock {\em Math. Res. Lett.}, 3(2):275--291, 1996.
\newblock Also \href{https://arxiv.org/abs/alg-geom/9506020}{alg-geom/9506020}.

\bibitem{gjh}
M.~Gross, D.~Huybrechts, and D.~Joyce.
\newblock {\em Calabi-{Y}au manifolds and related geometries}.
\newblock Universitext. Springer-Verlag, Berlin, 2003.
\newblock Lectures from the Summer School held in Nordfjordeid, June 2001.

\bibitem{egas}
A.~Grothendieck and J.~A. Dieudonn\'e.
\newblock {\em El\'ements de g\'eom\'etrie alg\'ebrique. {I}}, volume 166 of
  {\em Grundlehren der Mathematischen Wissen\-schaften}.
\newblock Springer-Verlag, Berlin, 1971.

\bibitem{dhl_magic_windows}
D.~Halpern-Leistner.
\newblock Derived {$\Theta$}-stratifications and the {$D$}-equivalence
  conjecture.
\newblock Preprint, \href{https://arxiv.org/abs/2010.01127}{arXiv:2010.01127}.

\bibitem{hartshorne}
R.~Hartshorne.
\newblock {\em Algebraic geometry}.
\newblock Springer-Verlag, New York-Heidelberg, 1977.
\newblock Graduate Texts in Mathematics, No. 52.

\bibitem{huybrechts_fm}
D.~Huybrechts.
\newblock {\em Fourier-{M}ukai transforms in algebraic geometry}.
\newblock Oxford Mathematical Monographs. The Clarendon Press Oxford University
  Press, Oxford, 2006.

\bibitem{huybrechts_k3}
D.~Huybrechts.
\newblock {\em Lectures on {K}3 surfaces}, volume 158 of {\em Cambridge Studies
  in Advanced Mathematics}.
\newblock Cambridge University Press, Cambridge, 2016.
\newblock Also available at
  \httpurl{www.math.uni-bonn.de/people/huybrech/K3Global.pdf}.

\bibitem{jl}
Q.~Jiang and N.~C. Leung.
\newblock Derived categories of projectivizations and flops.
\newblock {\em Adv. Math.}, 396:Paper No. 108169, 44, 2022.
\newblock Also \href{https://arxiv.org/abs/1811.12525}{arXiv:1811.12525}.

\bibitem{kapranov_gr}
M.~M. Kapranov.
\newblock Derived category of coherent sheaves on {G}rassmann manifolds.
\newblock {\em Izv. Akad. Nauk SSSR Ser. Mat.}, 48(1):192--202, 1984.

\bibitem{kawamata}
Y.~Kawamata.
\newblock {$D$}-equivalence and {$K$}-equivalence.
\newblock {\em J. Differential Geom.}, 61(1):147--171, 2002.
\newblock Also \href{https://arxiv.org/abs/math/0205287}{math/0205287}.

\bibitem{km}
J.~Koll\'{a}r and S.~Mori.
\newblock {\em Birational geometry of algebraic varieties}, volume 134 of {\em
  Cambridge Tracts in Mathematics}.
\newblock Cambridge University Press, Cambridge, 1998.

\bibitem{sandor}
S.~J. Kov\'{a}cs.
\newblock Singularities of stable varieties.
\newblock In {\em Handbook of moduli. {V}ol. {II}}, volume~25 of {\em Adv.
  Lect. Math. (ALM)}, pages 159--203. Int. Press, Somerville, MA, 2013.
\newblock Also \href{https://arxiv.org/abs/1102.1240}{arXiv:1102.1240}.

\bibitem{markman_brill_noether}
E.~Markman.
\newblock Brill-{N}oether duality for moduli spaces of sheaves on {K3}
  surfaces.
\newblock {\em J. Algebraic Geom.}, 10(4):623--694, 2001.
\newblock Also \href{https://arxiv.org/abs/math/9901072}{math/9901072}.

\bibitem{mukai_inventiones}
S.~Mukai.
\newblock Symplectic structure of the moduli space of sheaves on an abelian or
  {K3} surface.
\newblock {\em Invent. Math.}, 77(1):101--116, 1984.

\bibitem{mukai_tata}
S.~Mukai.
\newblock On the moduli space of bundles on {K3} surfaces. {I}.
\newblock In {\em Vector bundles on algebraic varieties ({B}ombay, 1984)},
  volume~11 of {\em Tata Inst. Fund. Res. Stud. Math.}, pages 341--413. Tata
  Inst. Fund. Res., Bombay, 1987.

\bibitem{nakajima}
H.~Nakajima.
\newblock Heisenberg algebra and {H}ilbert schemes of points on projective
  surfaces.
\newblock {\em Ann. of Math. (2)}, 145(2):379--388, 1997.
\newblock Also \href{https://arxiv.org/abs/alg-geom/9507012}{alg-geom/9507012}.

\bibitem{nakajima_k3}
H.~Nakajima.
\newblock Convolution on homology groups of moduli spaces of sheaves on {K3}
  surfaces.
\newblock In {\em Vector bundles and representation theory ({C}olumbia, {MO},
  2002)}, volume 322 of {\em Contemp. Math.}, pages 75--87. Amer. Math. Soc.,
  Providence, RI, 2003.

\bibitem{namikawa}
Y.~Namikawa.
\newblock Mukai flops and derived categories.
\newblock {\em J. reine angew. Math.}, 560:65--76, 2003.
\newblock Also \href{https://arxiv.org/abs/math/0203287}{math/0203287}.

\bibitem{negut}
A.~Negu\c{t}.
\newblock Hecke correspondences for smooth moduli spaces of sheaves.
\newblock {\em Publ. Math. Inst. Hautes \'{E}tudes Sci.}, to appear.
\newblock Also \href{https://arxiv.org/abs/1804.03645}{arXiv:1804.03645}.

\bibitem{og}
K.~O'Grady.
\newblock The weight-two {H}odge structure of moduli spaces of sheaves on a
  {K3} surface.
\newblock {\em J. Algebraic Geom.}, 6(4):599--644, 1997.
\newblock Also \href{https://arxiv.org/abs/alg-geom/9510001}{alg-geom/9510001}.

\bibitem{rrv}
J.~P. Ramis, G.~Ruget, and J.~L. Verdier.
\newblock Dualit\'{e} relative en g\'{e}om\'{e}trie analytique complexe.
\newblock {\em Invent. Math.}, 13:261--283, 1971.

\bibitem{stacks-project}
The {Stacks Project Authors}.
\newblock \textit{Stacks Project}.
\newblock \httpurl{stacks.math.columbia.edu}, 2020.

\bibitem{ryan_thesis}
R.~Takahashi.
\newblock {\em A categorical $\sl_2$ action on some moduli spaces of sheaves}.
\newblock PhD thesis, University of Oregon, 2020.

\bibitem{thomason}
R.~W. Thomason.
\newblock Les {K}-groupes d'un sch\'{e}ma \'{e}clat\'{e} et une formule
  d'intersection exc\'{e}dentaire.
\newblock {\em Invent. Math.}, 112(1):195--215, 1993.

\bibitem{yoshioka_reflection}
K.~Yoshioka.
\newblock Some examples of {M}ukai's reflections on {K3} surfaces.
\newblock {\em J. reine angew. Math.}, 515:97--123, 1999.

\bibitem{yoshioka_ade}
K.~Yoshioka.
\newblock An action of a {L}ie algebra on the homology groups of moduli spaces
  of stable sheaves.
\newblock In {\em Algebraic and arithmetic structures of moduli spaces
  ({S}apporo 2007)}, volume~58 of {\em Adv. Stud. Pure Math.}, pages 403--459.
  Math. Soc. Japan, Tokyo, 2010.

\end{thebibliography}
\bigskip
\scriptsize
\noindent Nicolas Addington \\
adding@uoregon.edu \\[\bigskipamount]
Ryan Takahashi \\
rtakahas@uoregon.edu \\[\bigskipamount]
Department of Mathematics \\
University of Oregon \\
Eugene, OR 97403-1222

\end{document}